\documentclass[preprint,12pt]{elsarticle}
\usepackage{hyperref}
\pdfoutput=1
\usepackage{enumerate}
\usepackage{graphics}
\usepackage{epsfig}
\usepackage{amsmath}
\usepackage{amsfonts}
\usepackage{amssymb}
\usepackage{latexsym}

\newtheorem{theorem}{Theorem}[section]
\newtheorem{lemma}{Lemma}[section]
\newtheorem{definition}{Definition}[section]

\newcommand{\ds}{\displaystyle}

\newcommand{\eps}{\varepsilon}

\usepackage{amssymb}

\begin{document}

\begin{frontmatter}


\title{ A parameter uniform essentially first order convergence of a fitted mesh method for a class of parabolic singularly perturbed Robin problem for 
a system of reaction-diffusion equations}

\author{R.Ishwariya}
\address{Department of Mathematics, Bishop Heber College, Tiruchirappalli, Tamil Nadu, India.}
\ead{ishrosey@gmail.com}

\author{J.J.H.Miller}
\address{Institute for Numerical Computation and Analysis, Dublin, Ireland.}
\ead{jm@incaireland.org}

\author{S.Valarmathi\corref{mycorrespondingauthor}}
\address{Department of Mathematics, Bishop Heber College, Tiruchirappalli, Tamil Nadu, India.}
\cortext[mycorrespondingauthor]{Corresponding author}
\ead{valarmathi07@gmail.com}

\begin{abstract}
In this paper, a class of linear parabolic systems of singularly perturbed second order differential equations of reaction-diffusion type with
initial and Robin boundary conditions is considered. The components of the solution $\vec u$ of this system exhibit parabolic boundary layers with sublayers. A numerical method composed of a classical finite difference scheme on a piecewise
uniform Shishkin mesh is suggested. This method is proved to be first order convergent in time and essentially first order convergent in the space variable in
the maximum norm uniformly in the perturbation parameters.
\end{abstract}

\end{frontmatter}
\section{Introduction}

\qquad A differential equation in which small parameters multiply the highest order derivative and some or none of the lower order derivatives
is known as a singularly perturbed differential equation. In this paper, a class of linear parabolic singularly perturbed second order differential equation of reaction-diffusion type with
initial and Robin boundary conditions is considered. \\

For a general introduction to parameter-uniform numerical methods for singular perturbation problems, see \cite{11}, \cite{12}, \cite{18} and \cite{19}. In \cite{13}, a Dirichlet boundary value problem for a linear parabolic singularly perturbed differential equation is studied and a numerical method 
comprising of a standard finite difference operator on a fitted piecewise uniform mesh is considered and it is proved to be uniform with respect to the small
parameter in the maximum norm. In \cite{14}, a boundary-value problem for a singularly perturbed parabolic PDE with convection is considered on an interval
in the case of the singularly perturbed Robin boundary condition is considered and using a defect correction technique, an ε-uniformly convergent schemes of
high-order time-accuracy is constructed. The efficiency of the new defect-correction schemes is confirmed by numerical experiments. In \cite{15}, a 
one-dimensional steady-state convection dominated convection-diffusion problem with Robin boundary conditions is considered and the numerical solutions 
obtained using an upwind finite difference scheme on Shishkin meshes are uniformly convergent with respect to the diffusion cofficient.\\

\indent Consider the following parabolic initial-boundary value problem for a singularly perturbed linear system of second order
differential equations
\begin{equation}
 \displaystyle\frac{\partial\vec{u}}{\partial t}(x,t)-E\frac{\partial^2\vec{u}}{\partial x^2}(x,t)+ A(x,t)\vec{u}(x,t)=\vec{f}(x,t), \text{ on } \Omega, \label{5p1}
 \end{equation}with
\begin{equation}\label{5p1aa} 
\begin{array}{l}
\vec{u}(0,t)-E_*\dfrac{\partial \vec u}{\partial x}(0,t)=\vec{\phi}_L(t),\;\;
 \vec{u}(1,t)+E_*\dfrac{\partial \vec u}{\partial x}(1,t)=\vec{\phi}_R(t),\;\;0\leq t \leq T,\\
\hspace{3.0cm} \vec{u}(x,0)=\vec{\phi}_B(x),\;\; 0\leq x \leq 1,
 \end{array}
 \end{equation}
where $\Omega = \{(x,t): 0<x<1,\; 0<t\leq T\},\;\; \bar{\Omega}=\Omega \cup \Gamma,\;\; \Gamma = \Gamma_L\cup\Gamma_B\cup\Gamma_R$ with
$\Gamma_L=\{(0,t) :0\leq t\leq T \},\;
\Gamma_R=\{(1,t) :0\leq t\leq T \}\text{ and } \Gamma_B=\{(x,0) :0 < x < 1 \}.$ Here, for all $(x,t)\in \bar
{\Omega},\;\vec{u}(x,t)\text{ and }\vec{f}(x,t)$ are column $n-$vectors, $E,$ $E_*$ and $A$ are $n \times n$ matrices,
$E = diag(\vec{\varepsilon}),\; \vec\varepsilon = (\varepsilon_1,...,\varepsilon_n),$
$E_* =  diag(\overrightarrow{\sqrt{\varepsilon}}),\; \overrightarrow{\sqrt{\varepsilon}} = (\sqrt\varepsilon_1,...,\sqrt\varepsilon_n)$ with $0 < \varepsilon_i < 1$ for all $i = 1,..., n.$ 
The parameters $\varepsilon_i$ are assumed to be distinct and for convenience, to have the ordering $\varepsilon_1 < ... <\varepsilon_n.$\\
The problem \eqref{5p1}, \eqref{5p1aa} can also be written in the operator form
\begin{equation*}L\vec{u} = \vec{f}\; \text{ on }\; \Omega,\end{equation*}
\begin{equation*}\beta_0\vec{u}(0,t)=\vec{\phi}_L(t),\;\;\beta_1\vec{u}(1,t)= \vec{\phi}_R(t),\;\;
\vec{u}(x,0)=\vec{\phi}_B(x),\end{equation*}
where the operators $L, \beta_0, \beta_1 $ are 
defined by \begin{equation*} L = \displaystyle I\frac{\partial}{\partial t}-E\frac{\partial^2}{\partial x^2}+ A,\;
\beta_0=I-E_*\dfrac{\partial}{\partial x},\;\beta_1=I+E_*\dfrac{\partial}{\partial x} \end{equation*}
where $I$ is the identity operator. The reduced problem corresponding to \eqref{5p1}, \eqref{5p1aa} is defined by
\begin{equation} \label{5p1b} \displaystyle\frac{\partial\vec{u}_0}{\partial t}+ A\vec{u}_0=\vec{f},\;  \text{ on } \;\Omega,\;\;
\vec{u}_0=\vec{u}\; \text{ on }\; \Gamma_B. \end{equation}
The problem \eqref{5p1}, \eqref{5p1aa} is said to be singularly perturbed in the following sense.\\
Each component $u_i,\;i=1,...,n$ of the solution $\vec u$ of \eqref{5p1}, \eqref{5p1aa} is expected to exhibit twin layers of width $O(\sqrt{\eps_n})$ at $x=0$ and $x=1$ 
while the components $u_i,\;i=1,...,n-1$ have additional twin layers of width $O(\sqrt{\eps_{n-1}}),$ the components $u_i,\;i=1,...,n-2$ have 
additional twin layers of width $O(\sqrt{\eps_{n-2}})$ and so on. 

\section{ Solution of the continuous problem}\label{c5s32}
\indent Standard theoretical results on the existence of the solution of \eqref{5p1}, \eqref{5p1aa} are stated, without proof, in this section. See \cite{16} and \cite{17} for more details.
For all $(x,t) \in \bar\Omega,$ it is assumed that the components $a_{ij} (x,t)$ of $A(x,t)$ satisfy the inequalities
\begin{equation} a_{ii} (x,t) >\displaystyle \sum^n_{^{j\neq i}_{j=1}}|a_{ij} (x,t)|\;\text{ for }\;1 \leq i \leq n,\;\text{ and }\;a_{ij} (x,t) \leq 0\;\text{ for}\;
i \neq j\label{5p3}\end{equation} and for some $\alpha,$  \begin{equation}\displaystyle0<\alpha<\min_{^{(x,t)\in \bar{\Omega}}_{1\leq i\leq n}}
(\sum^n_{j=1}a_{ij}(x,t)).\label{5p4}\end{equation}It is also assumed, without loss of generality, that
\begin{equation}\displaystyle\sqrt{\varepsilon_n}\leq\frac{\displaystyle\sqrt{\alpha}}{6}.\label{5p5} \end{equation}
\noindent Sufficient conditions for the existence, uniqueness and regularity of a solution of \eqref{5p1}, \eqref{5p1aa} are given in the following theorem.
\begin{theorem} \label{st} Assume that $ A \text{ and } \vec f$ are sufficiently smooth.  Also assume that ${\phi}_{L,i} \in C^{2}(\Gamma_L),\;{\phi}_{B,i} \in C^{5}(\Gamma_B)$,
\; ${\phi}_{R,i} \in C^{2}(\Gamma_R)$ and the following compatibility conditions are fulfilled at the corners $(0,0)$ and $(1,0)$ of\; $\Gamma.$\\
\begin{equation}\label{36} \vec{\phi}_B(0) = \vec{\phi}_L (0)+\dfrac{d \vec{\phi}_B}{dx}(0)\;\;\; \text{and}\;\;\; \vec{\phi}_B(1) = \vec{\phi}_R (0)-
\dfrac{d \vec{\phi}_B}{dx}(1),\end{equation}
\begin{equation}\label{37}
\begin{array}{lcl}
\ds\frac{d \vec{\phi}_L}{dt}(0)&=& -E \ds\frac{d^3 \vec{\phi}_B}{dx^3}(0)+E \ds\frac{d^2 \vec{\phi}_B}{dx^2}(0)+A(0,0) \ds\frac{d \vec{\phi}_
B}{dx}(0)-[A(0,0)-\dfrac{\partial A}{\partial x}(0,0)]\vec{\phi}_B(0)\\&&+\vec f (0,0)-\dfrac{\partial \vec f}{\partial x}(0,0),\\
\ds\frac{d \vec{\phi}_R}{dt}(0) &=& E \ds\frac{d^3 \vec{\phi}_B}{dx^3}(1)+E \ds\frac{d^2 \vec{\phi}_B}{dx^2}(1)-A(1,0) \ds\frac{d \vec{\phi}_
B}{dx}(1) -[A(1,0)+\dfrac{\partial A}{\partial x}(1,0)]\vec{\phi}_B(1)\\&&+\vec f (1,0)+\dfrac{\partial \vec f}{\partial x}(1,0),
\end{array}
\end{equation}
and
\begin{equation}\label{38}
\begin{array}{lcl}
\dfrac{d^2\vec{\phi}_L}{dt^2}(0)&=&-E^2\dfrac{d^5\vec{\phi}_B}{d x^5}(0)+E^2\dfrac{d^4\vec{\phi}_B}{dx^4}(0) + 2EA(0,0)\dfrac{d^3
\vec{\phi}_B}{dx^3}(0)+[-2E A(0,0)\\&&+4E\dfrac{\partial A}{\partial x}(0,0)]\dfrac{d^2\vec{\phi}_B}{dx^2}(0)
+[-2E\dfrac{\partial A}{\partial x}(0,0)+3E\dfrac{\partial^2 A}{\partial x^2}(0,0)-A^2(0,0)\\&&+\dfrac{\partial A}{\partial t}(0,0)]\dfrac{d\vec{\phi}_B}{dx}(0)
+[-E\dfrac{\partial^2 A}{\partial x^2}(0,0)+A^2(0,0)-\dfrac{\partial A}{\partial t}(0,0)+E\dfrac{\partial^3 A}{\partial x^3}(0,0)\\&&
-2A(0,0)\dfrac{\partial A}{\partial x}(0,0)+\dfrac{\partial^2 A}{\partial x \partial t}(0,0)]\vec{\phi}_B(0)+[-A(0,0)+\dfrac{\partial A}
{\partial x}(0,0)]\vec f(0,0)\\&&+\dfrac{\partial \vec f}{\partial t}(0,0)-E\dfrac{\partial^3 \vec f}{\partial x^3}(0,0)+E\dfrac{\partial^2 \vec f}{\partial x^2}(0,0)+A(0,0)\dfrac{\partial \vec f}{\partial x}(0,0)
-\dfrac{\partial^2 \vec f}{\partial x \partial t}(0,0),
\end{array}
\end{equation}
\begin{equation}\label{39}
\begin{array}{lcl}
\dfrac{d^2\vec{\phi}_R}{dt^2}(0)&=&E^2\dfrac{d^5\vec{\phi}_B}{d x^5}(1)+E^2\dfrac{d^4\vec{\phi}_B}{dx^4}(1) - 2EA(1,0)\dfrac{d^3
\vec{\phi_B}}{dx^3}(1)+[-2E A(1,0)\\&&-4E\dfrac{\partial A}{\partial x}(1,0)]\dfrac{d^2\vec{\phi}_B}{dx^2}(1)
+[-2E\dfrac{\partial A}{\partial x}(1,0)-3E\dfrac{\partial^2 A}{\partial x^2}(1,0)+A^2(1,0)\\&&-\dfrac{\partial A}{\partial t}(1,0)]\dfrac{d\vec{\phi}_B}{dx}(1)
+[-E\dfrac{\partial^2 A}{\partial x^2}(1,0)+A^2(1,0)-\dfrac{\partial A}{\partial t}(1,0)-E\dfrac{\partial^3 A}{\partial x^3}(1,0)\\&&
+2A(1,0)\dfrac{\partial A}{\partial x}(1,0)-\dfrac{\partial^2 A}{\partial x \partial t}(1,0)]\vec{\phi}_B(1)+[-A(1,0)-\dfrac{\partial A}
{\partial x}(1,0)]\vec f(1,0)\\&&+\dfrac{\partial \vec f}{\partial t}(1,0)+E\dfrac{\partial^3 \vec f}{\partial x^3}(1,0)+E\dfrac{\partial^2 \vec f}{\partial x^2}(1,0)
-A(1,0)\dfrac{\partial \vec f}{\partial x}(1,0)+\dfrac{\partial^2 \vec f}{\partial x \partial t}(1,0).
\end{array}
\end{equation}
\noindent Then there exists a unique solution $\vec u$ of \eqref{5p1}, \eqref{5p1aa} satisfying $ u_i \in C_\lambda ^{(4)} (\bar\Omega)$.
\end{theorem}
\section{Analytical results}\label{c5s33}
The operator $L$ satisfies the following maximum principle:
\begin{lemma}\label{5pmax}
Let the assumptions \eqref{5p3} - \eqref{5p5} hold. Let $\vec{\psi}$ be any vector-valued function in the domain of $L$ such that 
$\beta_0\vec\psi(0,t)\geq\vec0,\; \beta_1\vec\psi(1,t)\geq\vec0,\; \vec\psi(x,0)\geq\vec0.$ Then\; $L\vec\psi(x,t) \geq \vec0$\; on\; $\Omega$ \;implies that
$\vec\psi(x,t) \geq \vec0$\; on\; $\bar{\Omega}.$
\end{lemma}
\begin{lemma}\label{5psr} Let the assumptions \eqref{5p3} - \eqref{5p5} hold. If $\vec\psi$ is any vector-valued function in the domain of $L,$ then, \;for each\;
$i,\; 1 \leq i \leq n $ and\; $(x,t)\in\bar\Omega,$  $$|\psi_i(x,t)| \le\;
\max\displaystyle\left\{\parallel \beta_0\vec{\psi}(0,t)\parallel,\;\parallel \beta_1\vec{\psi}(1,t)\parallel,\;\parallel \vec{\psi}(x,0)\parallel,\; 
\dfrac{1}{\alpha}\parallel L\vec \psi\parallel\right\}.$$
\end{lemma}
A standard estimate of the solution $\vec u$ of the problem \eqref{5p1}, \eqref{5p1aa} and its derivatives is contained in the following lemma.
\begin{lemma}\label{5pud}Let the assumptions \eqref{5p3} - \eqref{5p5} hold and let $\vec u$ be the solution of \eqref{5p1}, \eqref{5p1aa}.
Then,\; for all $(x, t) \in \bar{\Omega}$ and each $i = 1, . . . , n,$\\
$|u_i(x,t)|\quad \leq C( \parallel \vec\phi_L(t) \parallel + \parallel \vec\phi_R(t) \parallel + \parallel \vec\phi_B(x) \parallel + \parallel \vec f \parallel ),\\
|\dfrac{\partial^l u_i}{\partial t^l}(x,t)| \leq C( \parallel \vec u \parallel +\ds\sum_{q=0}^{l} \parallel \dfrac{\partial^q\vec{f}}{\partial t^q} \parallel ),\;\; l= 1,2,\\
|\dfrac{\partial^l u_i}{\partial x^l}(x,t)| \leq C\eps_i^{\frac{-l}{2}}( \parallel \vec{u} \parallel + \parallel \vec{f} \parallel + \parallel \dfrac{\partial\vec{f}}{\partial t} \parallel ), \;\; l=1,2,\\
|\dfrac{\partial^l u_i}{\partial x^l}(x,t)| \leq
C\eps^{-1}_i \eps^{\frac{-(l-2)}{2}}_1 ( \parallel \vec{u} \parallel + \parallel \vec{f} \parallel +
 \parallel \dfrac{\partial\vec{f}}{\partial t} \parallel + \parallel \dfrac{\partial^2 \vec f}{\partial t^2} \parallel +
\eps^{\frac{l-2}{2}}_1 \parallel \dfrac{\partial^{l-2} \vec f}{\partial x^{l-2}} \parallel ), l=3,4,\\
|\dfrac{\partial^l u_i}{\partial x^{l-1} \partial t}(x,t)| \leq C\eps_i^{\frac{-(l-1)}{2}}( \parallel \vec{u} \parallel + \parallel \vec{f} \parallel + \parallel \dfrac{\partial\vec{f}}{\partial
t} \parallel + \parallel \dfrac{\partial^2 \vec f}{\partial t^2} \parallel ), \;\; l=2,3.$
\end{lemma}
\noindent The Shishkin decomposition of the solution $\vec u$ of the problem \eqref{5p1}, \eqref{5p1aa} is 
\begin{equation}\label{5p15a}\vec u=\vec v+\vec w \end{equation}
where $\vec v$ and $\vec w$ are the smooth and singular components of the solution $\vec u$ respectively.\\
Taking into consideration, the sublayers that appear for the components, the smooth component $\vec v$ is subjected to further decomposition. 
\begin{equation}\label{5p15aa}
\begin{array}{lcl}
 v_n&=&u_{0,n}+\eps_n v_{n,n},\\
 v_{n-1}&=&u_{0,{n-1}}+\eps_n v_{n-1,n}^1,\\
 \vdots\\
 v_1&=&u_{0,1}+\eps_n v_{1,n}^1,
\end{array}
 \end{equation}
as all the components have $\eps_n$ layers. Since components except $u_n$ have $\eps_{n-1}$ sublayers, the components $v_{n-1},...,v_1$ takes the form,
\begin{equation}\label{5p15ab}
\begin{array}{lcl}
v_{n-1}&=&u_{0,{n-1}}+\eps_n (v_{n-1,n}+\eps_{n-1}v_{n-1,n-1}),\\
 v_{n-2}&=&u_{0,{n-2}}+\eps_n (v_{n-2,n}+\eps_{n-1}v_{n-2,n-1}^1),\\
 \vdots\\
 v_1&=&u_{0,1}+\eps_n (v_{1,n}+\eps_{n-1}v_{1,n-1}^1).
\end{array}
 \end{equation}
Further, $u_{n-2},\;u_{n-3},...,u_2,u_1$ have $\eps_{n-2}$ sublayers and hence that leads to the decomposition,
\begin{equation}\label{5p15ac}
\begin{array}{rcl}
 v_{n-2}&=&u_{0,{n-2}}+\eps_n (v_{n-2,n}+\eps_{n-1}(v_{n-2,n-1}+\eps_{n-2}v_{n-2,n-2})),\\
 v_{n-3}&=&u_{0,{n-3}}+\eps_n (v_{n-3,n}+\eps_{n-1}(v_{n-3,n-1}+\eps_{n-2}v_{n-3,n-2}^1)),\\
\vdots\\
 v_1&=&u_{0,1}+\eps_{n}(v_{1,n}+\eps_{n-1}(v_{1,n-1}+\eps_{n-2}v_{1,n-2}^1)).
\end{array}
 \end{equation}
Proceeding like this, it is not hard to see that
\begin{equation*}\label{5p16aa}
\begin{pmatrix}
v_1\\
v_2\\
 \vdots\\
v_n
\end{pmatrix}
=\begin{pmatrix}
u_{0,1}\\
u_{0,2}\\
 \vdots\\
u_{0,n}
\end{pmatrix}
\begin{pmatrix}
 \gamma_1\\
 \gamma_2\\
 \vdots\\
\gamma_n
\end{pmatrix}
\end{equation*}i.e.
\begin{equation}\label{5p16ab}
 \vec v(x,t)=\vec{u}_0(x,t)+\vec\gamma(x,t)
\end{equation}where
\begin{equation}\label{5p16ac}
\gamma_j= \vec{\eps\;}^n(\vec {v\;}_j^j)^T,
\end{equation}
\[\vec{\eps\,}^n=(\eps_1\eps_{2}...\eps_n,\;\eps_{2}\eps_{3}...\eps_n,......,\;\eps_{n-1}\eps_{n},\;\eps_{n}),\;\;
\vec{v\,}^i_i=(0,\;0,...,\;v_{i,i},\;v_{i,i+1},......,\;v_{i,n}).\]
Then using \eqref{5p15a} and \eqref{5p16ab} in \eqref{5p1}, \eqref{5p1aa}, it is found that the smooth component $\vec v$
of the solution $\vec u$ satisfies
\begin{equation}\label{5p15b}
L\vec{v}=\vec{f}, \text{ on } \Omega 
 \end{equation}with
\begin{equation}\label{5p15c} 
\qquad\qquad\beta_0\vec{v}(0,t)=\beta_0(\vec{u}_0+\vec{\gamma})(0,t),\;\;\beta_1\vec{v}(1,t)= \beta_1(\vec{u}_0+\vec{\gamma})(1,t),\;\;\vec v(x,0)=(\vec u_0+\vec\gamma)(x,0).\end{equation}
and the singular component $\vec w$ of the solution $\vec u$ satisfies
 \begin{equation}
 L\vec{w}=\vec{0}, \text{ on } \Omega \label{5p16g}
 \end{equation}with
\begin{equation}\label{5p16h} 
\begin{array}{rcl}
\beta_0\vec{w}(0,t)=\beta_0(\vec{u}-\vec{v})(0,t),\;\;
\beta_1\vec{w}(1,t)=\beta_1(\vec{u}-\vec{v})(1,t),\;\;
\vec w(x,0)=\vec 0.
 \end{array}
 \end{equation}
 Consider the following parabolic initial-boundary value problem for a singularly perturbed linear system of second order
differential equations
\begin{equation}
\displaystyle\frac{\partial\vec{\hat u}}{\partial t}(x,t)-\hat{E}\frac{\partial^2\vec{\hat u}}{\partial x^2}(x,t)+ \hat A(x,t)\vec{\hat u}(x,t)=\vec{\hat f}(x,t), \text{ on } 
\Omega, \label{5pn1}
 \end{equation}with
\begin{equation}\label{5pn1aa} \begin{array}{l}
\hat{u}_2(0,t)-\sqrt{\eps_n}\dfrac{\partial \hat{u}_n}{\partial x}(0,t)=\alpha(t),\;
\hat{u}_2(1,t)+\sqrt{\eps_n}\dfrac{\partial \hat{u}_n}{\partial x}(1,t)=\beta(t),\;0\leq t \leq T,\\
 \hspace{3.0cm}\vec{\hat u}(x,0)=\vec{\delta}(x),\; 0\leq x \leq 1,\;
 \end{array}
\end{equation}
where $\hat E$ is a $n \times n$ matrix,
$\hat E = diag(0,0,...,0,\varepsilon_n)$ with $0 < \varepsilon_n < 1.$\\
The problem \eqref{5pn1}, \eqref{5pn1aa} can also be written in the operator form
\begin{equation*}\hat L\vec{\hat u} = \vec{\hat f}\; \text{ on }\; \Omega, \label{4hn1m}\end{equation*}
\begin{equation*}\begin{array}{c}b_0\hat{u}_n(0,t)=\alpha(t),\;\;b_1\hat{u}_n(1,t)= \beta(t),\;\;
\vec{\hat u}(x,0)=\vec{\delta}(x).\label{5pn2m}\end{array}\end{equation*}
where the operators $\hat L, b_0, b_1 $ are 
defined by \begin{equation*} \hat L = \displaystyle I\frac{\partial}{\partial t}-\hat E\frac{\partial^2}{\partial x^2}+ \hat A,\;
b_0=I-\sqrt{\eps_n}\dfrac{\partial}{\partial x},\;b_1=I+\sqrt{\eps_n}\dfrac{\partial}{\partial x}  \end{equation*}
where $I$ is the identity operator. The reduced problem corresponding to \eqref{5pn1}, \eqref{5pn1aa} is defined by
\begin{equation*} \label{4hn1b} \displaystyle\frac{\partial\vec{\hat u}_0}{\partial t}+ \hat A\vec{\hat u}_0=\vec{\hat f},\;  \text{ on } \;\Omega,\;\;
\vec{\hat u}_0=\vec{\hat u}\; \text{ on }\; \Gamma_B. \end{equation*}
 The operator $\hat L$ satisfies the following maximum principle:
\begin{lemma}\label{3qnmax}
 Let the assumptions \eqref{5p3} - \eqref{5p5} hold. Let $\vec{\psi}=(\psi_1,...,\psi_n)^T$ be any vector-valued function in the domain of $\hat L$ such that 
$b_0\psi_n(0,t)\geq 0,\; b_1\psi_n(1,t)\geq 0,\; \vec\psi(x,0)\geq\vec0.$ Then\; $\hat L\vec\psi(x,t) \geq \vec0$\; on\; $\Omega$ \;implies that
$\vec\psi(x,t) \geq \vec0$\; on\; $\bar{\Omega}.$
\end{lemma}
\begin{lemma}\label{5pnsr} Let the assumptions \eqref{5p3} - \eqref{5p5} hold. If $\vec\psi=(\psi_1,...,\psi_n)^T$
is any vector-valued function in the domain of $\hat L,$ then, \;for each $i=1,...,n$ and $(x,t)\in\bar\Omega,$
\begin{equation*}\begin{array}{c}
|\psi_i(x,t)| \le\;\max\displaystyle\left\{\parallel b_0{\psi}(0,t)\parallel,\parallel b_1{\psi}(1,t)\parallel,
\parallel \vec{\psi}(x,0)\parallel, 
\dfrac{1}{\alpha}\parallel \hat L\vec \psi\parallel\right\}.
\end{array}
\end{equation*}
\end{lemma}
A standard estimate of the solution $\vec{\hat u}$ of the problem \eqref{5pn1}, \eqref{5pn1aa} and its derivatives is contained in the following lemma.
\begin{lemma}\label{5pnud}Let the assumptions \eqref{5p3} - \eqref{5p5} hold and let $\vec{\hat u}$ be the solution of \eqref{5pn1}, \eqref{5pn1aa}.
Then,\; for all $(x, t) \in \bar{\Omega}$ and each $i = 1,...,n,$\\
$|\hat u_i(x,t)|\quad \leq C( \parallel \alpha(t) \parallel + \parallel \beta(t) \parallel + \parallel \vec\delta(x) \parallel + \parallel \vec{\hat f} \parallel ),\\
|\dfrac{\partial^l \hat u_i}{\partial t^l}(x,t)| \leq C( \parallel \vec{\hat u} \parallel +\sum_{q=0}^{l} \parallel \dfrac{\partial^q\vec{\hat f}}{\partial t^q} \parallel ),\;\; l= 1,2,\\
|\dfrac{\partial \hat u_i}{\partial x}(x,t)| \leq C\eps_n^{-1/2} ( \parallel \vec{\hat u} \parallel + \parallel \vec{\hat f} \parallel +\parallel \dfrac{\partial\vec{\hat f}}{\partial t} \parallel + \eps_n^{1/2}\parallel \dfrac{\partial\vec{\hat f}}{\partial x} \parallel ),\\
|\dfrac{\partial^2 \hat u_i}{\partial x^2}(x,t)| \leq C\eps_n^{-1}( \parallel \vec{\hat u} \parallel + \parallel \vec{\hat f} \parallel +\parallel \dfrac{\partial\vec{\hat f}}{\partial t} \parallel + \eps_n\parallel \dfrac{\partial\vec{\hat f}}{\partial x} \parallel+ \eps_n\parallel \dfrac{\partial^2\vec{\hat f}}{\partial x^2} \parallel ),\\
|\dfrac{\partial^3 \hat u_i}{\partial x^3}(x,t)| \leq C\eps_n^{-3/2}( \parallel \vec{\hat u} \parallel + \parallel \vec{\hat f} \parallel +\parallel \dfrac{\partial\vec{\hat f}}{\partial t} \parallel+\parallel \dfrac{\partial^2\vec{\hat f}}{\partial t^2} \parallel + \eps_n^{3/2}\parallel \dfrac{\partial\vec{\hat f}}{\partial x} \parallel+ \eps_n^{3/2}\parallel \dfrac{\partial^2\vec{\hat f}}{\partial x^2} \parallel )+ \eps_n^{3/2}\parallel \dfrac{\partial^3\vec{\hat f}}{\partial x^3} \parallel ),\\
|\dfrac{\partial^4 \hat u_i}{\partial x^4}(x,t)| \leq C \eps_n^{-2} ( \parallel \vec{\hat u} \parallel + \parallel \vec{\hat f} \parallel +\parallel \dfrac{\partial\vec{\hat f}}{\partial t} \parallel+\parallel \dfrac{\partial^2\vec{\hat f}}{\partial t^2} \parallel + \eps_n^{2}\parallel \dfrac{\partial\vec{\hat f}}{\partial x} \parallel+ \eps_n^{2}\parallel \dfrac{\partial^2\vec{\hat f}}{\partial x^2} \parallel )+ \eps_n^{2}\parallel \dfrac{\partial^3\vec{\hat f}}{\partial x^3} \parallel )+ \eps_n^{2}\parallel \dfrac{\partial^4\vec{\hat f}}{\partial x^4} \parallel ),\\
\|\dfrac{\partial^2 \hat u_i}{\partial x \partial t}(x,t)| \leq C\eps_n^{-1/2}( \parallel \vec{\hat u} \parallel + \parallel \vec{\hat f} \parallel + \parallel \dfrac{\partial \vec{\hat f}}{\partial t} \parallel + \parallel \dfrac{\partial^2 \vec{\hat f}}{\partial t^2} \parallel )+ \eps_n^{1/2}\parallel \dfrac{\partial\vec{\hat f}}{\partial x} \parallel ),\\
|\dfrac{\partial^3 \hat u_i}{\partial x^2 \partial t}(x,t)| \leq C\eps_n^{-1}( \parallel \vec{\hat u} \parallel + \parallel \vec{\hat f} \parallel + \parallel \dfrac{\partial \vec{\hat f}}{\partial t} \parallel + \parallel \dfrac{\partial^2 \vec{\hat f}}{\partial t^2} \parallel )+ \eps_n\parallel \dfrac{\partial\vec{\hat f}}{\partial x} \parallel+ \eps_n\parallel \dfrac{\partial^2\vec{\hat f}}{\partial x^2} \parallel ).$
\end{lemma}
\noindent Bounds on the smooth component $\vec{v}$ of $\vec{u}$ and its derivatives are contained in
\begin{lemma}\label{5pvd} Let the assumptions \eqref{5p3} - \eqref{5p5} hold. Then there exists a constant $C,$ such that, for each $(x,t) \in \bar \Omega$ 
and $i=1,..., n$,
\begin{equation*}\label{5p16k} \begin{array} {lcl}
|\dfrac{\partial^l v_i}{\partial t^l}(x,t)|\leq C, \;\; l=0,1,2,\;\;
 &&|\dfrac{\partial^l v_i}{\partial x^l}(x,t)|\leq C, \;\; l=1,2,\\
 |\dfrac{\partial^l v_i}{\partial x^l}(x,t)|\leq C \eps_{i}^{-(l-2)/2},\;\; l=3,4,\;\;
&&|\dfrac{\partial^{l+1} v_i}{\partial x^l\partial t}(x,t)|\leq C,\;\; l=1,2.
\end{array} \end{equation*}
\end{lemma}
\noindent \textbf{Proof.} 
From \eqref{5p15aa} - \eqref{5p15ac} it is observed that the components $v_{i,j},\;\;i=1,...,n,\;\;j=i,i+1,...,n$ satisfy the following systems of equations:
\begin{equation}\label{5p17a}
\begin{array}{rcl}
\dfrac{\partial v_{1,n}}{\partial t}+a_{11}v_{1,n}+a_{12}v_{2,n}+...+a_{1n}v_{n,n}&=&\dfrac{\eps_1}{\eps_n}\dfrac{\partial^2 u_{0,{1}}}{\partial x^2}\\
\dfrac{\partial v_{2,n}}{\partial t}+a_{21}v_{1,n}+a_{22}v_{2,n}+...+a_{2n}v_{n,n}&=&\dfrac{\eps_2}{\eps_n}\dfrac{\partial^2 u_{0,{2}}}{\partial x^2}\\
\vdots&&\\
\dfrac{\partial v_{n-1,n}}{\partial t}+a_{n-11}v_{1,n}+a_{n-12}v_{2,n}+...+a_{n-1n}v_{n,n}&=&\dfrac{\eps_{n-1}}{\eps_n}\dfrac{\partial^2 u_{0,{n-1}}}{\partial x^2}\\
\dfrac{\partial v_{n,n}}{\partial t}-\eps_n\dfrac{\partial^2 v_{n,n}}{\partial x^2}+a_{n1}v_{1,n}+a_{n2}v_{2,n}+...+a_{nn}v_{n,n}&=&\dfrac{\partial^2 u_{0,{n}}}{\partial x^2}
\end{array}
\end{equation}
with
\begin{equation}\label{5p17b}
 (v_{n,n}-\sqrt{\eps_n}\dfrac{\partial v_{n,n}}{\partial x})(0,t)=0,\;\;(v_{n,n}+\sqrt{\eps_n}\dfrac{\partial v_{n,n}}{\partial x})(1,t)=0,\;\; v_{i,n}(x,0)=0,
 \end{equation}
where $u_{0,i},\;i=1,...,n$ is the solution of the reduced problem \eqref{5p1b}.
 \begin{equation}\label{5p17c}
\begin{array}{rcl}
\dfrac{\partial v_{1,n-1}}{\partial t}+a_{11}v_{1,n-1}+...+a_{1n-1}v_{n-1,n-1}&=&\dfrac{\eps_1}{\eps_{n-1}}\dfrac{\partial^2v_{1,n}}{\partial x^2}\\ 
\dfrac{\partial v_{2,n-1}}{\partial t}+a_{21}v_{1,n-1}+...+a_{2n-1}v_{n-1,n-1}&=&\dfrac{\eps_2}{\eps_{n-1}}\dfrac{\partial^2v_{2,n}}{\partial x^2}\\
\vdots&&\\
\dfrac{\partial v_{n-2,n-1}}{\partial t}+a_{n-21}v_{1,n-1}+...+a_{n-2n-1}v_{n-1,n-1}&=&\dfrac{\eps_{n-2}}{\eps_{n-1}}\dfrac{\partial^2v_{n-2,n}}{\partial x^2}\\
\dfrac{\partial v_{n-1,n-1}}{\partial t}-\eps_{n-1}\dfrac{\partial^2v_{n-1,n-1}}{\partial x^2}+a_{n-11}v_{1,n-1}+...&+&a_{n-1n-1}v_{n-1,n-1}\\
&=&\dfrac{\partial^2v_{n-1,n}}{\partial x^2}
\end{array}
\end{equation}
with
\begin{equation}\label{5p17d}
\begin{array}{rcl}
(v_{n-1,n-1}-\sqrt{\eps_{n-1}}\dfrac{\partial v_{n-1,n-1}}{\partial x})(0,t)=0,\;
&&(v_{n-1,n-1}+\sqrt{\eps_{n-1}}\dfrac{\partial v_{n-1,n-1}}{\partial x})(1,t)=0,\\
 \qquad\qquad v_{i,n-1}(x,0)=0,
\end{array}
\end{equation} and so on.\\
Lastly,
\begin{equation}\label{5p17e}
\begin{array}{rcl}
\dfrac{\partial v_{1,2}}{\partial t}+a_{11}v_{1,2}+a_{12}v_{2,2}&=&\dfrac{\eps_1}{\eps_2}\dfrac{\partial^2v_{1,3}}{\partial x^2}\\
\dfrac{\partial v_{2,2}}{\partial t}-\eps_2 \dfrac{\partial^2v_{2,2}}{\partial x^2}+a_{21}v_{1,2}+a_{22}v_{2,2}&=&\dfrac{\partial^2v_{2,3}}{\partial x^2}
\end{array}
\end{equation}
with
\begin{equation}\label{5p17f}
 (v_{2,2}-\sqrt{\eps_{2}}\dfrac{\partial v_{2,2}}{\partial x})(0,t)=0,\;\;(v_{2,2}+\sqrt{\eps_{2}}\dfrac{\partial v_{2,2}}{\partial x})(1,t)=0,\;\; v_{i,2}(x,0)=0,
 \end{equation} and 
 \begin{equation}\label{5p17g}
\dfrac{\partial v_{1,1}}{\partial t}- \eps_1\dfrac{\partial^2 v_{1,1}}{\partial x^2}+a_{11}v_{1,1}=\dfrac{\partial^2v_{1,2}}{\partial x^2}
 \end{equation} with
 \begin{equation}\label{5p17h}
 (v_{1,1}-\sqrt{\eps_{1}}\dfrac{\partial v_{1,1}}{\partial x})(0,t)=0,\;\;(v_{1,1}+\sqrt{\eps_{1}}\dfrac{\partial v_{1,1}}{\partial x})(1,t)=0,\;\; v_{1,1}(x,0)=0.
 \end{equation}
 From the expressions \eqref{5p17a}-\eqref{5p17h} and using Lemma \eqref{5pnud} for $\vec v$, it is found that for $ i=1,...,n,\;\;j=i,i+1,...,n,\;\;i\leq j,\;\;k=1,2,3,4,\;\;l=0,1,2,\;\;m=1,2$
\begin{equation}\label{5p17i}
\begin{array}{lcl}
 |\dfrac{\partial^l v_{i,j}}{\partial t^l}(x,t)|\leq C(1+\ds\prod_{r=j+1}^{n}\eps_r^{-1}),\;\;\;\;
 |\dfrac{\partial^k v_{i,j}}{\partial x^k}(x,t)|\leq C(1+\eps_{j}^{-k/2}\ds\prod_{r=j+1}^{n}\eps_r^{-1}),\\
\hspace{3.0cm}|\dfrac{\partial^{m+1} v_{i,j}}{\partial x^{m}\partial t}(x,t)|\leq C(1+\eps_{j}^{-m/2}\ds\prod_{r=j+1}^{n}\eps_r^{-1}).
\end{array}
\end{equation}
From \eqref{5p16ab}, \eqref{5p16ac} and \eqref{5p17i}, the following bounds for $v_i,\;\;i=1,2,...,n$ hold:
\begin{equation*}
\begin{array}{rcl}
|\dfrac{\partial^l v_i}{\partial t^l}(x,t)|\leq C, \;\; l=0,1,2,\;\;
 &&|\dfrac{\partial^l v_i}{\partial x^l}(x,t)|\leq C, \;\; l=1,2,\\
 |\dfrac{\partial^l v_i}{\partial x^l}(x,t)|\leq C \eps_{i}^{-(l-2)/2},\;\; l=3,4,\;\;
&&|\dfrac{\partial^{l+1} v_i}{\partial x^l\partial t}(x,t)|\leq C,\;\; l=1,2.
\end{array}
\end{equation*}\\
\noindent The layer functions $B^{L}_{i},\; B^{R}_{i},\; B_{i},\; i=1,\;
\dots,\; n,$ associated with the solution $\vec u,$ are
defined on $\bar\Omega$ by
\[B^{L}_{i}(x) = e^{-x\sqrt{\alpha/\varepsilon_i}},\;B^{R}_{i}(x) =
B^{L}_{i}(1-x),\;B_{i}(x) = B^{L}_{i}(x)+B^{R}_{i}(x).\] 
The following elementary properties of these layer functions, for all $0 \leq x <
y \leq 1,$ should be noted:\\
$B_i(x)=B_i(1-x),
B^{L}_1(x)< B^{L}_2(x),\;B^{L}_1(x)> B^{L}_2(y), \;0<B^{L}_i(x)\leq1,
B^{R}_1(x) < B^{R}_2(x),\;B^{R}_1(x) <B^{R}_2(y), \;0<B^{R}_i(x)\leq1,
B_{i}(x)\; \text{is monotonically decreasing for increasing}\; x \in [0,\frac{1}{2}],\\
B_{i}(x) \text{\;is monotonically increasing for increasing}\; x \in [\frac{1}{2},1],
B_{i}(x) \leq 2B_{i}^L(x)\; \text{for}\; x \in [0,\frac{1}{2}],\\
 \;B_{i}(x) \leq 2B_{i}^R(x) \;\text{for}\; x \in [\frac{1}{2},1],
B^L_i(2\frac{\sqrt\varepsilon_i}{\sqrt\alpha}\ln N)=N^{-2}.$\\
The interesting points $x^{(s)}_{i,j}$ are now defined.
\begin{definition}\label{5pdefn1.1}
For $B_i^L$, $B_j^L,$ each $i,j, \;1 \leq i \neq j \leq n$ and
each $s,\; s>0,$ the point $x^{(s)}_{i,j}$ is defined by
\begin{equation*}\label{5p31}
\frac{B^L_i(x^{(s)}_{i,j})}{\varepsilon^s _i}=
\frac{B^L_j(x^{(s)}_{i,j})}{\varepsilon^s _j}. \end{equation*}
It is remarked that
\begin{equation*}\label{5p32}
\frac{B^R_i(1-x^{(s)}_{i,j})}{\varepsilon^s_i}=
\frac{B^R_j(1-x^{(s)}_{i,j})}{\varepsilon^s_j}. \end{equation*} \end{definition}

In the next lemma, the existence, uniqueness and ordering of the points $x^{(s)}_{i,j}$ are established. 
Sufficient conditions for them to lie in the domain $\bar\Omega$ are also provided.
\begin{lemma}\label{5plem1.6} For all $\,i,j\,$ such that $1 \leq i < j \leq
n$ and $0 < s \leq 3/2,$ the points $x_{i,j}^{(s)}$ exist, are uniquely
defined and satisfy the following inequalities
\begin{equation*}\label{5p33}
\frac{B^L_{i}(x)}{\varepsilon^s_i} > \frac{B^L_{j}(x)}{\varepsilon^s_j},\;\; x \in [0,x^{(s)}_{i,j}),\;\; \frac{B^L_{i}(x)}{\varepsilon^s_i} <
\frac{B^L_{j}(x)}{\varepsilon^s_j}, \; x \in (x^{(s)}_{i,j}, 1].\end{equation*}
In addition, the following ordering holds 
\begin{equation*}\label{5p34}
x^{(s)}_{i,j} < x^{(s)}_{i+1,j}, \; \text{if} \;\; i+1<j \;\;
\text{and} \;\; x^{(s)}_{i,j} <
x^{(s)}_{i,j+1}, \;\; \text{if} \;\; i<j. \end{equation*} 
Also,
\begin{equation*}\label{5p35}
 x^{(s)}_{i,j}< 2s\frac{\sqrt\varepsilon_j}{\sqrt\alpha}\;\; \text{and} \;\;
x^{(s)}_{i,j} \in (0,\frac{1}{2})\;\; \text{if} \;\; i<j. \end{equation*}
Analogous results hold for $B^R_i,\; B^R_j$ and the points $1-x^{(s)}_{i,j}.$
\end{lemma}
\noindent \textbf{Proof.} The proof is as given in \cite{21}.\\
\noindent Bounds on the singular component $\vec{w}$ of $\vec{u}$ and its derivatives are contained in
\begin{lemma}\label{5pwd} Let the assumptions \eqref{5p3} - \eqref{5p5} hold. Then there exists a constant $C,$ such that, for each $(x,t) \in \bar{\Omega} $ and $i=1,..., n$,
\begin{equation*}\begin{array} {l}
|\dfrac{\partial^l w_i}{\partial t^l}(x,t)| \le C B_{n}(x),\;\text{ for} \;\; l=0,1,2, \;\;
|\dfrac{\partial^l w_i}{\partial x^l}(x,t)| \le C\ds\sum_{r=i}^n \dfrac{B_{r}(x)}{\eps_r^{\frac{l}{2}}},\;\text{ for}\;\; l=1,2,\\
|\dfrac{\partial^3 w_i}{\partial x^3}(x,t)| \le C\ds\sum_{r=1}^n \dfrac{B_{r}(x)}{\eps_r^{\frac{3}{2}}},\;\;
|\eps_i\dfrac{\partial^4 w_i}{\partial x^4}(x,t)| \le C \ds\sum_{r=1}^n \dfrac{B_{r}(x)}{\eps_r}.
\end{array} \end{equation*}
\end{lemma}
\noindent \textbf{Proof.} 
To derive the bound of $\vec{w}$, define $\vec{\psi}^\pm(x,t)=(\psi_1,...,\psi_n)^T,$ where
$${\psi_i}^\pm(x,t)=Ce^{\alpha t}B_n(x)\;\pm\;w_i(x,t), \text{ for each } i=1,\dots,n.$$
For a proper choice of $C,$\; $\beta_0\vec{\psi}^\pm(0,t)\geq \vec 0$,\;$\beta_1\vec{\psi}^\pm(1,t)\geq \vec 0$\;and\;${\vec{\psi}^\pm}(x,0)\geq \vec 0.$ Also, for\;$(x,t)\in \Omega, \;
L \vec \psi^\pm(x,t)\ge \vec 0.$ By Lemma \ref{5pmax}, $\vec\psi^{\pm} \geq \vec0$ on $\bar\Omega$ and it follows that
\begin{equation*}\label{5p1414}
 |w_i(x,t)|\leq Ce^{\alpha t}B_n(x)\;\;\;\text{or}\;\;\; |w_i(x,t)|\leq CB_n(x).
\end{equation*}
Differentiating the homogeneous equation satisfied by $w_i$, partially with respect to
$`t`$, and using Lemma \ref{5pmax}, it is not hard to see that
 \begin{equation*}\label{5p1416}
|\frac{\partial w_i}{\partial t}(x,t)|\leq CB_n(x).
\end{equation*}
Note that,\[|\frac{\partial^2 w_i}{\partial x \partial t}(x,t)|\leq |\frac{\partial^2 u_i}{\partial x \partial t}(x,t)|+|\frac{\partial^2 v_i}{\partial x \partial t}(x,t)|.\]
Thus,\[ |\frac{\partial^2 w_i}{\partial x \partial t}(x,t)| \leq C {\eps_i}^{\frac{-1}{2}}({ \parallel \vec u \parallel }+ \parallel \vec f \parallel + \parallel  \frac{\partial \vec f}
{\partial t} \parallel + \parallel \frac{\partial^2 \vec f}{\partial t^2} \parallel ).\]
Similarly,\begin{equation*}\label{5p1417}  |\frac{\partial^3 w_i}{\partial x^2 \partial t}(x,t)| \leq C \eps_i^{-1}({ \parallel \vec u \parallel }+ \parallel \vec f \parallel + \parallel  \frac{\partial \vec f}{\partial t} \parallel + \parallel \frac{\partial^2 \vec f}{\partial t^2} \parallel ). \end{equation*}
As before, using suitable barrier functions, it is not hard to verify that
\[ |\frac{\partial^{l+1} w_i}{\partial x^l \partial t}(x,t)| \leq C \eps_i^{\frac{-l}{2}}B_n(x),\;\;l=1,2.\]
Differentiating the equation satified by $w_i$ partially with respect to $`t`$ once and rearranging, yields
\begin{equation*}
|\frac{\partial^2 w_i}{\partial t^2}(x,t)|\leq CB_n(x).
\end{equation*}
\noindent The bounds on $\dfrac{\partial^l w_i}{\partial x^l},\; l=1,2,3,4$ and
$i=1,\dots,n$ are now derived by induction on $n$. For $n=1$, the result follows from . It is then
assumed that the required bounds on $\dfrac{\partial w_i}{\partial x} ,
\dfrac{\partial^2 w_i}{\partial x^2},\dfrac{\partial^3
w_i}{\partial x^3}$\;and\;$\dfrac{\partial^4 w_i}{\partial x^4}$
hold for all systems up to order $n-1$. Define
$\vec{\tilde{w}}=(w_1,\dots,w_{n-1})$, then
$\vec{\tilde{w}}$
 satisfies the system 
 \begin{equation}\label{5p4444}
 \frac{\partial \vec{\tilde{w}}}{\partial t}-\tilde{E}\frac{\partial^2 \vec{\tilde{w}}}{\partial x^2}+\tilde A \vec{\tilde{w}} = \vec g, 
\end{equation}
with
\begin{equation*}\label{5p4444a} 
\begin{array}{lcl}
\beta_0\vec{\tilde{w}}(0,t) = \beta_0(\vec{\tilde{u}}-\vec{\tilde{v}})(0,t),\;\;
\beta_1\vec{\tilde{w}}(1,t) = \beta_1(\vec{\tilde{u}}-\vec{\tilde{v}})(1,t),\;\;
\vec{\tilde{w}}(x,0) = \vec{\tilde{0}}.
\end{array}
\end{equation*}
Here, $\tilde{E}$ and $\tilde{A}$ are the matrices obtained by
deleting the last row and last column from $E,A$ respectively, the
components of $\vec g$ are $g_i = -a_{in}w_n$ for $1\leq i \leq
n-1$ and $\vec{\tilde{v}}=\vec{\tilde{u}}_0+\vec{\tilde{\gamma}}$ is the corresponding component decomposition of $\vec{\tilde v}$ similar to \eqref{5p16ab} of $\vec{v}.$
Now decompose $\vec{\tilde{w}}$ into smooth and singular
components to get $ \vec{\tilde{w}} = \vec{p} +\vec{q},$ where
$L \vec{p} = \vec{g},\;\; \beta_0\vec p(0,t)=\beta_0(\vec {\tilde{u}}_0+\vec {\tilde{\gamma}})(0,t), \;\;\beta_1\vec p(1,t)=\beta_1(\vec {\tilde{u}}_0+\vec {\tilde{\gamma}})(1,t),\;\;
\vec p(x,0) = (\vec {\tilde{u}}_0+\vec {\tilde{\gamma}})(x,0)\;\text{ and }\;L\vec q=\vec 0,\;\; \beta_0\vec q(0,t) = \beta_0\vec {\tilde{w}}(0,t) - \beta_0\vec p(0,t),  \;\;
\beta_1\vec q(1,t) = \beta_1\vec {\tilde{w}}(1,t) - \beta_1\vec p(1,t),\;\;\vec q(x,0) = \vec{\tilde{w}}(x,0)-\vec {p}(x,0).\\$ 
Consider the equation of the system satisfied by $w_i,$
\[\frac{\partial w_i}{\partial t} - \eps_i \frac{\partial^2 w_i}{\partial x^2 } + \sum_{j=1}^{n} a_{ij}{ w_j}=0. \]
By using mean-value theorem, the bound on $\dfrac{\partial w_i}{\partial x },$ for each $(x,t)$ is 
determined as follows:\begin{equation*} \label{5p1422}
|\frac{\partial w_i}{\partial x}(x,t)|\leq C \eps_i^{\frac{-1}{2}}B_n(x).
\end{equation*}
Rearranging the equation of the system satisfied by $w_i,$ yields
\begin{equation*}\label{5p1418a}|\dfrac{\partial^2 w_i}{\partial x^2}(x,t)|\leq C \eps_i^{-1} B_i(x).\end{equation*}
Differentiating the equation satisfied by $w_i$ with respect to $`x`$ once and twice and rearranging, the following bounds are derived
\[|\frac{\partial^3 w_i}{\partial x^3}(x,t)| \leq C \sum_{r=1}^{n} {\eps_r}^{-\frac{3}{2}}B_r(x), \;\;
|\eps_i\frac{\partial^4 w_i}{\partial x^4}(x,t)| \leq C \sum_{r=1}^{n} {\eps_r}^{-1}B_r(x).\]
\indent Using the bounds on $w_n, \dfrac{\partial w_n}{\partial x},\dfrac{\partial^2 w_n}{\partial x^2},\dfrac{\partial^3 w_n}{\partial x^3}$ and $\dfrac{\partial^4
w_n}{\partial x^4}$, it is seen that the function $\vec g$ in \eqref{5p4444} and its derivatives $\dfrac{\partial \vec g}{\partial x},$ $\dfrac{\partial^2 \vec g}{\partial x^2},$ $\dfrac{\partial^3 \vec g}{\partial
x^3},$ $\dfrac{\partial^4 \vec g}{\partial x^4}$ are bounded by $ CB_n(x),$ $C\dfrac{B_n(x)}{\sqrt \eps_n},$ $C\dfrac{B_n(x)}{\eps_n},$ $C\sum_{r=1}^{n}
\dfrac{B_r(x)}{\eps_r^{\frac{3}{2}}},$ and $C \eps_n^{-1} \sum_{r=1}^{n}\dfrac{B_r(x)}{\eps_r} $ respectively.
Introducing the functions ${\vec \psi}^\pm(x,t)=C e^{\alpha t}B_n(x)\vec e\;\pm\;\vec p(x,t)$, it is easy to see that ${\beta_0\vec {\psi}^\pm}(0,t)=Ce^{\alpha t}B_n(0)\vec e\;\pm\;\beta_0\vec p(0,t) \geq
\vec0$, ${\beta_1\vec {\psi}^\pm}(1,t)=Ce^{\alpha t}B_n(1)\vec e\;\pm\;\beta_1\vec p(1,t) \geq \vec 0$,  ${\vec \psi}^\pm(x,0)=CB_n(x)\vec e\;\pm\;\vec p(x,0) \geq \vec 0$ and
\[(L {\vec \psi}^\pm)_i(x,t)=C(-\eps_i \frac{\alpha}{\eps_n}+\alpha e^{\alpha t}+\sum_{j=1}^{n}a_{ij})B_n(x)\;\pm\;(L\vec p)_i \geq 0,\; \text{as} -\frac{\eps_i}{\eps_n}\geq -1.\]
Applying Lemma \ref{5pmax}, it follows that $\parallel\vec p(x,t)\parallel \leq C B_n(x)$.\\ Defining the barrier functions through ${\vec \theta}^\pm (x,t) = C \eps_n^{-\frac{l}{2}}e^{\alpha t}B_n(x) \vec
e\;\pm\; \dfrac{\partial^l \vec p}{\partial x^l},\;l=1,2 $ and using Lemma \ref{5pmax} for the problem satisfied by $\vec p$ and the bounds of the derivatives of $\vec{g}$, the bounds of $\dfrac{\partial \vec p}{\partial x} $
and $\dfrac{\partial^2 \vec p}{\partial x^2} $ are derived.\\ The bounds for $ \dfrac{\partial^l \vec p}{\partial x^l},\;l=3,4$ follow from the defining equation of $\vec p.$\\
\noindent By induction, the following bounds for $\vec q$ hold for $i=1,\dots,n-1,$
\begin{equation*}
\begin{array}{lclclcl}
 |\dfrac{\partial q_i}{\partial x}(x,t)| &\leq& C\ds\left[\frac{B_i(x)}{{\sqrt \eps_i}}+\dots+\frac{B_{n-1}(x)}{{\sqrt \eps_{n-1}}}\right],\\
 |\dfrac{\partial^2 q_i}{\partial x^2} (x,t)|&\leq& C\ds\left[\frac{B_i(x)}{{ \eps_i}}+\dots+\frac{B_{n-1}(x)}{{ \eps_{n-1}}}\right],\\
 |\dfrac{\partial^3 q_i}{\partial x^3} (x,t)|&\leq& C\ds\left[\frac{B_1(x)}{{ \eps_1^{3/2}}}+\dots+\frac{B_{n-1}(x)}{{ \eps_{n-1}^{3/2}}}\right],\\
 |\eps_i\dfrac{\partial^4 q_i}{\partial x^4}(x,t) |&\leq& C\ds\left[\frac{B_1}{{ \eps_1}}+\dots+\frac{B_{n-1}(x)}{{ \eps_{n-1}}}\right].
\end{array}
\end{equation*}
Combining the bounds for the derivatives of $p_i$ and $q_i,$ it follows that, for $i= 1,2,\dots,n-1,$
\[ \begin{array}{l}
|\dfrac{\partial^l \tilde{w}_i}{\partial x^l}(x,t)|\leq C \sum_{r=i}^{n-1} \dfrac{B_r(x)}{{ \eps_r^{\frac{l}{2}}}}\;\;\text{for}\;\; l= 1,2,\\
|\dfrac{\partial^3 \tilde{w}_i}{\partial x^3}(x,t)|\leq C \sum_{r=1}^{n-1} \dfrac{B_r(x)}{ \eps_r^{\frac{3}{2}}},\;\;
|\eps_i\dfrac{\partial^4 \tilde{w}_i}{\partial x^4}(x,t)|\leq C \sum_{r=1}^{n-1} \dfrac{B_r(x)}{\eps_r}. \end{array}\] 
 
 Using the above bounds along with the bounds of $w_i$ and its derivatives, the proof of the lemma for the system of $n$
 equations gets completed.
\section{The Shishkin mesh}\label{c5s34}
A piecewise uniform Shishkin mesh is now
constructed. Let $\Omega^M_t=\{t_k \}_{k=1}^{M},\;\Omega^N_x=\{x_j
\}_{j=1}^{N-1},\;\bar\Omega^M_t=\{t_k
\}_{k=0}^{M},\;\bar\Omega^N_x=\{x_j
\}_{j=0}^{N},\;\Omega^{N,M}=  \Omega^N_x \times \Omega^M_t,\;
 \bar{\Omega}^{N,M}=\bar\Omega^N_x \times \bar\Omega^M_t \;\;\text{ and }\;\;\Gamma^{N,M}=\Gamma \cap \bar{\Omega}^{N,M}.$
 The mesh $\bar\Omega^M_t$ is chosen to be a uniform mesh with $M$
sub-intervals on $[0,T]$. The mesh $\bar\Omega^N_x$ is a
piecewise-uniform Shishkin mesh with $N$ mesh intervals. The interval $[0,1]$ is subdivided into 
$5$ sub-intervals given by
\[[0,\sigma_1]\cup(\sigma_{1},\sigma_2]\cup(\sigma_2,1-\sigma_2]\cup(1-\sigma_2,1-\sigma_{1}]\cup(1-\sigma_1,1].\]
The parameters $\sigma_1$ and $\sigma_2$ which determine the points separating
the uniform meshes, are defined by 
\begin{equation}\label{3151}\sigma_{2} = 
\min\displaystyle\left\{\frac{1}{4},\; 2\frac{\sqrt\varepsilon_2}{\sqrt\alpha}\ln
N\right\}\;\;\text{and}\;\;\sigma_{1}=\min\displaystyle\left\{\frac{\sigma_{2}}{2},
\; 2\frac{\sqrt\varepsilon_1}{\sqrt\alpha}\ln
N\right\}.\end{equation} Also, $\sigma_0=0,\;\sigma_{3}=\frac{1}{2}.$
Clearly $0<\sigma_1<\sigma_2\le\frac{1}{4}, \qquad
\frac{3}{4}\leq 1-\sigma_2 < 1-\sigma_1 <1.$
Then, on the sub-interval $(\sigma_2,1-\sigma_2]$ a uniform mesh with
$\frac{N}{2}$ mesh-points is placed and on each of the
sub-intervals $[0,\sigma_1],\;[\sigma_1,\sigma_{2}),\;[1-\sigma_{2},1-\sigma_1)\;\text{and}\;[1-\sigma_1,1)$
a uniform mesh of $\frac{N}{8}$ mesh-points is placed. \\
Thus $\bar\Omega^{N,M}$ is a piecewise uniform Shishkin grid with $NM$ mesh elements.\\
In practice, it is convenient to take $ N = 8 q,\;q 
\geq 3.$\\
In particular, when the parameters $\sigma_1,\;\sigma_2,$ are with the left choice, the Shishkin mesh 
$\bar{\Omega}^{N,M}$ becomes the classical uniform mesh with the transition parameters $\sigma_1 = \dfrac{1}
{8},\;\sigma_2 = \dfrac{1}{4}$ and with the step size $N^{-1}$ throughout on $\bar{\Omega}_x^N$.
\\
The Shishkin mesh suggested here has the following features: (i) when all the transition parameters have the left choice, it is the classical uniform mesh and (ii) it is coarse in the outer region and becomes finer and finer towards the left and right boundaries. From the above construction it is clear that the transition points
$\{\sigma_1, \sigma_2, 1-\sigma_1, 1-\sigma_2\}$ on $\bar{\Omega}^N_x$ are the only points at which the mesh-size can change and that it does not necessarily change at each of these points.
The following notations are introduced: if $x_j=\sigma_r$, then $h_r^-=x_j-x_{j-1},\; h_r^+=x_{j+1}-x_j$, $J= \{x_j=\sigma_r, 1-\sigma_r: h_r^+ \neq  h_r^-\}$.
In general, for each point $x_j$ in the sub-interval $(\sigma_{1}
,\sigma_2 ]$ and $(1-\sigma_{2},1-\sigma_1 ],$
\begin{equation}\label{3153} x_j-x_{j-1}
=8N^{-1}(\sigma_2-\sigma_1).\end{equation} Also, for
$x_j \in(\sigma_2,1-\sigma_2], \; x_j-x_{j-1}=2N^{-1}(1-2\sigma_2)$ and for
$x_j \in (0,\sigma_1]$ and $x_j\in (1-\sigma_1,1) , \;  x_j-x_{j-1}=8N^{-1}\sigma_1.$ Thus, for
 $r=1,2,$ the
change in the mesh-size at the point $x_j=\sigma_r$ is $h^+_r-h^-_r=8N^{-1}(d_r-d_{r-1}),$ where $d_r =\dfrac{r \sigma_{r+1}}{r+1}-\sigma_r$
with the convention
$d_0 =0.$ Notice that $d_r \ge 0$, $\bar\Omega^{N,M}$ is a
classical uniform mesh when $d_r = 0$ for $ r=1,2$
and, 
 from \eqref{3151} that
\begin{equation}\label{3156}\sigma_r \leq C \sqrt{\eps_r} \ln N,
 \; \;\; r=1,2. \end{equation}
It follows from (\ref{3153}) and (\ref{3156}) that for $r = 1,$
\begin{equation*}\label{3157}
h_{r}^-+h_r^+ \leq C\sqrt{\eps_{r+1}}N^{-1}\ln N.\end{equation*}
Also,
\begin{equation*}\label{3158}
\sigma_1=\frac{\sigma_{2}}{2},\;\text{ when}\; d_1=d_2
=0.
\end{equation*}
\begin{lemma}\label{5ps1} Assume that $d_r>0$ for some $r, 1 \leq r \leq n.$  Then the following
inequalities hold
\begin{equation*}\label{5p159}
B^L_r(1-\sigma_r) \leq B^L_r(\sigma_r)=N^{-2};\;\;
x^{(s)}_{r-1,r}\;\leq\;\sigma_r-h_r^- \;\text{ for} \; 
0 < s \leq 3/2.\end{equation*}
\begin{equation*}\label{5p1511} B_q^L(\sigma_r-h_r^-)\leq
CB_q^L(\sigma_r)\;\;\text{ for}\;\; 1 \leq r \leq q \leq n;\;\;
\frac{B_{q}^{L}(\sigma_{r})}{\sqrt{\eps_q}}\leq C\frac{1}{\sqrt{\eps_r}\ln N} \;\;
\text{ for} \;\; 1 \leq q \leq n.
\end{equation*}
Analogous results hold for $B_r^R$.
\end{lemma}
\noindent \textbf{Proof.} The proof is as given in \cite{21}.

\section{The discrete problem}\label{c5s35}
In this section a classical finite difference operator with an
appropriate Shishkin mesh is used to construct a numerical method
for the problem \eqref{5p1}, \eqref{5p1aa} which is shown later to be first order parameter-uniform convergent in time and essentially first order
parameter-uniform convergent in the space variable.

The discrete initial-boundary value problem is now defined by the finite difference scheme on the Shishkin mesh $\bar\Omega^{N,M},$ defined in the previous section.
\begin{equation}\label{5p161} 
 D^-_t\vec{U}(x_j,t_k)-E\delta^2_x\vec{U}(x_j,t_k) +A(x_j,t_k)\vec{U}(x_j,t_k)=\vec{f}(x_j,t_k)\;\; \text{ on }\;\; \Omega^{N,M},
 \end{equation}with
 \begin{equation}\label{5p161a}
 \begin{array}{lcl}
 \vec U(0,t_k) - {E_*}D^+_x \vec U(0,t_k) = \vec \phi_L(t_k),\;\;
 \vec U(1,t_k) + {E_*}D^-_x \vec U(1,t_k) = \vec \phi_R(t_k)\\
\hspace{4.0cm} \vec{U}(x_j,0) = \vec \phi_B(x_j).
\end{array}
 \end{equation}
The problem \eqref{5p161}, \eqref{5p161a} can also be written in the operator form
\[L^{N,M} \vec{U}=\vec{f} \;\; \text{ on }
\;\; \Omega^{N,M},\]
\begin{equation*}
 \begin{array}{lcl}
 \beta_0^{N,M}\vec{U}(0,t_k) = \vec \phi_L(t_k)\;\;
 \beta_1^{N,M}\vec{U}(1,t_k) = \vec \phi_R(t_k)\;\;
 \vec{U}(x_j,0) = \vec \phi_B(x_j),
\end{array}
 \end{equation*}
where\[L^{N,M} = I D^-_t-E\delta^2_x+A,\; \beta_0^{N,M}= I -E_*D^+_x,\;\beta_1^{N,M}= I +E_*D^-_x\]
The following discrete results are analogues to those for the
continuous case.\\
\begin{lemma}\label{5pdmax} Let the assumptions \eqref{5p3} - \eqref{5p5} hold.
Then, for any vector-valued mesh function $\vec{\Psi}$, the inequalities 
$\beta_0^{N,M}\vec{\Psi}(0,t_k) \ge \vec 0,\;\beta_1^{N,M}\vec{\Psi}(1,t_k) \ge \vec 0,\;\vec{\Psi}(x_j,0) \ge \vec 0$  and $L^{N,M}
\vec{\Psi}\;\ge\;\vec 0$ on\; $\Omega^{N,M}$ imply that $\vec
{\Psi}\ge \vec 0$ on $\bar{\Omega}^{N,M}.$
\end{lemma}
\noindent An immediate consequence of this is the following discrete stability
result.
\begin{lemma}\label{5pdsr} Let the assumptions \eqref{5p3} - \eqref{5p5} hold.
Then, for any vector-valued mesh function $\vec{\Psi}$ defined on $\bar\Omega^{N,M}$ and $i=1, \dots ,n $,
\[|{\Psi}_i(x_j,t_k)|\;\le\;\max\left\{ \parallel \beta_0^{N,M}\vec{\Psi}(0,t_k) \parallel ,\; \parallel \beta_1^{N,M}\vec{\Psi}(1,t_k) \parallel ,\; \parallel \vec{\Psi}(x_j,0) \parallel ,\; \frac{1}{\alpha} \parallel 
L^{N,M}\vec{\Psi} \parallel \right\}. \]
\end{lemma}
\noindent The following comparison principle will be used in the proof of the error estimate.
\begin{lemma}\label{5pcomp}Assume that,
for the vector-valued mesh functions $\vec{\Phi}$ and  $\vec{Z}$ satisfy\;\;
$|\beta_0^{N,M}\vec Z(0,t_k)| \leq \beta_0^{N,M}\vec\Phi(0,t_k),\;
|\beta_1^{N,M}\vec Z(1,t_k)| \leq \beta_1^{N,M}\vec\Phi(1,t_k),\;| \vec Z(x_j,0)| 
\leq \vec\Phi(x_j,0),\text{ and }|L^{N,M}\vec{Z}| \leq L^{N,M} \vec{\Phi}\;\text{ on }\;\Omega^{N,M}.$ Then, 
$|\vec Z| \leq \vec\Phi\;\text{ on }\; \bar\Omega^{N,M}.$
\end{lemma}
\section{The local truncation error}\label{c5s36}
From Lemma \ref{5pdsr}, it is seen that in order to bound the error
$\vec{U}-\vec{u}$, it suffices to bound
$\beta_0^{N,M}(\vec{U}-\vec u)(0,t_k),\;\beta_1^{N,M}(\vec{U}-\vec u)(1,t_k),\;(\vec U-\vec u)(x_j,0) \text{ and }L^{N,M}(\vec{U}-\vec{u})$. Note that, for $(x_j, t_k) \in \Omega ^{N,M},$
$L^{N,M}(\vec{U}-\vec{u})=L^{N,M}\vec{U}-L^{N,M}\vec{u} = \vec{f}-L^{N,M}\vec{u}
=L\vec{u}-L^{N,M}\vec{u}
=(L-L^{N,M})\vec{u}.$\\
It follows that, $ L^{N,M}(\vec{U}-\vec{u})=(\frac{\partial}{\partial
t}-D^-_t)\vec{u}-E(\frac{\partial^2}{\partial
x^2}-\delta^2_x)\vec{u}.$\\
Let $\vec{V}, \vec{W}$ be the discrete analogues of
$\vec{v}, \vec{w}$ respectively, given by
\begin{equation*}
L^{N,M}\vec V = \vec f \text{ on } \Omega^{N,M},
\end{equation*}
\begin{equation*}
\beta_0^{N,M}\vec{V}(0,t_k)=\beta_0\vec{v}(0,t_k),\; \beta_1^{N,M}\vec{V}(1,t_k)=\beta_1\vec{v}(1,t_k),\; \vec{V}(x_j,0)=\vec{v}(x_j,0), 
\end{equation*}
\begin{equation*}
L^{N,M}\vec{W} = \vec 0 \text{ on } \Omega^{N,M},
\end{equation*}
\begin{equation*}
\beta_0^{N,M}\vec{W}(0,t_k)=\beta_0\vec{w}(0,t_k),\;  \beta_1^{N,M}\vec{W}(1,t_k)=\beta_1\vec{w}(1,t_k),\; \vec{W}(x_j,0)=\vec{w}(x_j,0),
\end{equation*}
where $\vec v$ and $\vec w$ are the solutions of \eqref{5p15b}, \eqref{5p15c} and \eqref{5p16g}, \eqref{5p16h} respectively.\\
Therefore, the local truncation error of the smooth and singular components can be treated separately. Note that, for any smooth function $ \psi $ and for each
$(x_j,t_k)\in \Omega^{N,M}$, the following
distinct estimates of the local truncation error hold:
\begin{equation}\label{3174g}
|(\frac{\partial}{\partial t}-D^-_t)\psi(x_j,t_k)|\;\le\;
C(t_k-t_{k-1})\max_{s\;\in\;[t_{k-1},\;t_k]}|\frac{\partial^2\psi}{\partial
t^2}(x_j,s)|,
\end{equation}
\begin{equation}\label{3175ag}
|(\frac{\partial}{\partial x}-D^-_x)\psi(x_j,t_k)|\;\le\;
C(x_{j}-x_{j-1})\max_{s\;\in\;[x_{j-1},\;x_{j}]}|\frac{\partial^2\psi}{\partial
x^2}(s, t_k)|,
\end{equation}
\begin{equation}\label{3175g}
|(\frac{\partial}{\partial x}-D^+_x)\psi(x_j,t_k)|\;\le\;
C(x_{j+1}-x_{j})\max_{s\;\in\;[x_{j},\;x_{j+1}]}|\frac{\partial^2\psi}{\partial
x^2}(s, t_k)|,
\end{equation}
\begin{equation}\label{3176g}
\hspace{-2.5cm}|(\frac{\partial^2}{\partial x^2}-\delta^2_x)\psi(x_j,t_k)|\;\le\;
C\max_{s\;\in\;I_j}|\frac{\partial^2\psi}{\partial x^2}(s,t_k)|,
\end{equation}
\begin{equation}\label{3177g}
\hspace{-.5cm}|(\frac{\partial^2}{\partial
x^2}-\delta^2_x)\psi(x_j,t_k)|\;\le\;C(x_{j+1}-x_{j-1}) \max_{s\;\in\;I_j}|\frac{\partial^3\psi}{\partial x^3}(s,t_k)|.
\end{equation}
Here $I_j=[x_{j-1}, x_{j+1}]$.
\section{Error estimate}\label{c5s37}
The proof of the theorem on the error estimate is broken into two parts. First, a theorem concerning the error in the smooth component is established.
Then the error in the singular component is estimated.\\
\noindent 
Define the barrier function through
\begin{equation*}\vec{\Phi}(x_j,t_k)=C[M^{-1}+(r+1)N^{-1}\ln N+
(N^{-1}\ln N)\ds\sum_{\{r: \;\sigma_r \in
J\}}\dfrac{\sigma_r}{\sqrt{\eps_i}}\theta_r (x_j,t_k)]\vec{e},\end{equation*} where $C$ is sufficiently large and 
$\theta_r$ is a piecewise linear polynomial for each $x_j=\sigma_r\in J$ defined by
\begin{equation*} \theta_r(x,t)=
\left\{ \begin{array}{l}\;\; \dfrac{x}{\sigma_r}, \;\; 0 \leq x \leq \sigma_r,\\
\;\; 1, \;\; \sigma_r < x < 1-\sigma_r,\\ 
\;\; \dfrac{1-x}{\sigma_r}, \;\; 1-\sigma_r \leq x \leq 1.  \end{array}\right .\end{equation*}
 Also note that,\begin{equation}\label{5p181}(L^{N,M}\theta_r\vec{e})_i(x_j,t_k) \ge
\left\{ \begin{array}{l}\;\;
\alpha \theta_{r}(x_j, t_k), \qquad \quad \;\;  \text{ if} \;\;x_j \notin J \\
\;\; \alpha+\dfrac{2\eps_i}{ \sigma_r
(h_r^-+h_r^+)}, \;\;  \text{ if} \;\;x_j \in J.\end{array}\right .\end{equation}
Then, on $\Omega^{N,M}$, the components $\Phi_i$ of  $\vec{\Phi}$ satisfy
\begin{equation*}\label{5p183}0 \leq \Phi_{i}(x_j,t_k) \leq
C(M^{-1}+N^{-1}\ln N),\;\; 1 \leq i \leq n.\end{equation*}
Also,
\begin{equation}\label{5p184a}
( \beta_0\vec\Phi)_i(0,t) \ge C(M^{-1}+N^{-1}\ln N),\;\;
(\beta_1\vec\Phi)_i(1,t) \ge C(M^{-1}+N^{-1}\ln N)
\end{equation}
For $x_j \notin J$, it is not hard to see that,
\begin{equation}\label{5p184}
(L^{N,M}\vec{\Phi})_i(x_j,t_k) \ge C(M^{-1}+N^{-1}\ln
N)\end{equation} and, for $x_j \in J,$ it is not hard to see that,
\begin{equation}\label{5p185}
(L^{N,M}\vec{\Phi})_i(x_j,t_k) \ge C(M^{-1}+N^{-1}\ln N).\end{equation}
The following theorem gives the estimate of the error in the smooth component $\vec V$.
  \begin{theorem}\label{5psmootherrorthm} Let the assumptions \eqref{5p3} - \eqref{5p5} hold. Let $\vec v$ denote the smooth component of the 
solution of the problem \eqref{5p1}, \eqref{5p1aa} and $\vec V$ denote the smooth component of the
solution of the problem \eqref{5p161}, \eqref{5p161a}.  Then
\begin{equation*}\;\;  \parallel \vec{V}-\vec{v} \parallel  \leq C(M^{-1}+N^{-1}\ln N). \end{equation*}
\end{theorem}
\noindent In order to estimate the error in the singular component $\vec W$, the following lemmas are required.\\
\begin{lemma} \label{5pwwe}Assume that $x_j \notin J$.  Let the assumptions \eqref{5p3} - \eqref{5p5} hold. Then,
on $\Omega^{N,M}$, for each $1 \leq i \leq n$,
\begin{equation*}\label{5p1813} |(L^{N,M}(\vec{W}-\vec {w}))_i(x_j,t_k)|\leq C(M^{-1}+\frac{(x_{j+1}-x_{j-1})}{\sqrt{\eps_1}}).\end{equation*}
\end{lemma}
\noindent The following decomposition in the singular components $w_i$ are used in the next lemma
\begin{equation}\label{5p1814}w_i=\sum_{m=1}^{r+1}w_{i,m},\end{equation} where the components $w_{i,m}$ 
are defined by
\[\hspace{-3cm}w_{i,r+1}=\left\{ \begin{array}{ll} p^{(s)}_{i} & \text{on}\;\;[0,x^{(s)}_{r,r+1})\\
w_i & \text{on}\;\;[x^{(s)}_{r,r+1}, 1-x^{(s)}_{r,r+1}]\\
q^{(s)}_{i} & \text{on}\;\;(1-x^{(s)}_{r,r+1},1]\end{array}\right. \]
\noindent where  \[\hspace{-1.2cm} p^{(s)}_{i}(x,t)=\left\{ \begin{array}{ll}\sum_{k=0}^3
\dfrac{\partial^k w_i}{\partial
x^k}(x^{(s)}_{r,r+1},t)\dfrac{(x-x_{r,r+1}^{(s)})^k}{k!},\;\text{for}\;s=\frac{3}{2},\\
\sum_{k=0}^4
\dfrac{\partial^k w_i}{\partial
x^k}(x^{(s)}_{r,r+1},t)\dfrac{(x-x_{r,r+1}^{(s)})^k}{k!},\;\text{for}\;s=1,\end{array}\right.\]
\[ q^{(s)}_{i}(x,t)=\left\{ \begin{array}{ll}\sum_{k=0}^3
\dfrac{\partial^k w_i}{\partial
x^k}(1-x^{(s)}_{r,r+1},t)\dfrac{(x-(1-x_{r,r+1}^{(s)}))^k}{k!},\;\text{for}\;s=\frac{3}{2},\\ 
\sum_{k=0}^4
\dfrac{\partial^k w_i}{\partial
x^k}(1-x^{(s)}_{r,r+1},t)\dfrac{(x-(1-x_{r,r+1}^{(s)}))^k}{k!},\;\text{for}\;s=1,\end{array}\right.\]

\noindent and, for each $m$,  $r \ge m \ge 2$,
\[\hspace{-1cm}w_{i,m}=\left\{ \begin{array}{ll} p^{(s)}_{i} & \text{on} \;\; [0,x^{(s)}_{m-1,m})\\
w_i-\ds\sum_{k=m+1}^{r+1} w_{i,k} & \text{on}\;\;[x^{(s)}_{m-1,m}, 1-x^{(s)}_{m-1,m}]\\
q^{(s)}_{i} & \text{on} \;\; (1-x^{(s)}_{m-1,m},1]
\end{array}\right. \]

\noindent where  \[\hspace{-1.2cm} p^{(s)}_{i}(x,t)=\left\{ \begin{array}{ll}\sum_{k=0}^3
\dfrac{\partial^k w_i}{\partial
x^k}(x^{(s)}_{m,m+1},t)\dfrac{(x-x_{m,m+1}^{(s)})^k}{k!},\;\text{for}\;s=\frac{3}{2},\\ 
\sum_{k=0}^4
\dfrac{\partial^k w_i}{\partial
x^k}(x^{(s)}_{m,m+1},t)\dfrac{(x-x_{m,m+1}^{(s)})^k}{k!},\;\text{for}\;s=1,\end{array}\right.\]
\[ q^{(s)}_{i}(x,t)=\left\{ \begin{array}{ll}\sum_{k=0}^3
\dfrac{\partial^k w_i}{\partial
x^k}(1-x^{(s)}_{m,m+1},t)\dfrac{(x-(1-x_{m,m+1}^{(s)}))^k}{k!},\;\text{for}\;s=\frac{3}{2},\\ 
\sum_{k=0}^4
\dfrac{\partial^k w_i}{\partial
x^k}(1-x^{(s)}_{m,m+1},t)\dfrac{(x-(1-x_{m,m+1}^{(s)}))^k}{k!},\;\text{for}\;s=1,\end{array}\right.\]
 and \[ w_{i,1}=w_i-\sum_{k=2}^{r+1} w_{i,k}\;\; \text{on} \;\; [0,1]. \]
\noindent Notice that the decomposition \eqref{5p1814} depends on the choice of the polynomials $p^{(s)}_{i},\;q^{(s)}_{i}$ and the definition of $x^{(s)}_{i,j},\;1-x^{(s)}_{i,j}$ 
given by \eqref{5p31} and \eqref{5p32}.\\
\noindent The following lemma provides estimates of the derivatives of the components $w_{i,m},\; 1\leq m \leq r+1 $ of $w_i, \; 1\leq i\leq n.$
\begin{lemma} \label{5pdwd}Assume that $d_r>0$ for some $r$, $1 \leq r \leq n.$  
Let the assumptions \eqref{5p3} - \eqref{5p5} hold. Then, for
each $1 \leq i \leq n$,   the components in the decomposition \eqref{5p1814}
satisfy the following
estimates for each $q$ and $r$,  $1 \le q \le r$, and all
$(x_j, t_k) \in \Omega^{N,M},$
\begin{eqnarray*}
\begin{array}{l}
|\dfrac{\partial^2 w_{i,q}}{\partial x^2}(x_j,t_k)|  \leq C \min \{\dfrac{1}{\eps_q}, \dfrac{1}{\eps_i}\}B_q(x_j),\;\;
|\dfrac{\partial^3 w_{i,q}}{\partial x^3}(x_j,t_k)|  \leq C \min \{ \dfrac{1}{\eps_i \sqrt{\eps_q}},\dfrac{1}{\eps_q^{3/2}}\}B_q(x_j),\\
|\dfrac{\partial^3 w_{i,r+1}}{\partial x^3}(x_j,t_k)|   \leq C \min \{ \sum_{q=r+1}^n \dfrac{B_q(x_j)}{\eps_i \sqrt{\eps_q}}, \sum_{q=r+1}^n \dfrac{B_q(x_j)}{\eps_q^{3/2}}\},\;\;
|\dfrac{\partial^4 w_{i,q}}{\partial x^4}(x_j,t_k)|   \leq C \dfrac{B_q(x_j)}{\eps_i \eps_q},\\
|\dfrac{\partial^4 w_{i,r+1}}{\partial x^4}(x_j,t_k)|   \leq C\sum_{q=r+1}^n \dfrac{B_q(x_j)}{\eps_i \eps_q}.
\end{array}
\end{eqnarray*}
\end{lemma}
\begin{lemma}\label{5pdwd1} Assume that $d_r>0$ for some $r$, $1 \leq r \leq n$. Let the assumptions \eqref{5p3} - \eqref{5p5} hold.
Then,\\ if $x_j \notin J,$
\begin{equation}\label{5p1815}
|(L^{N,M}(\vec{W}-\vec{w}))_i(x_j,t_k)| \leq C[M^{-1}+B_r(x_{j-1})+\dfrac{(x_{j+1}-x_{j-1})}{\sqrt{\eps_{r+1}}}],
\end{equation}
and if $x_j \in J,$
\begin{equation} \label{5p1816}
|(L^{N,M}(\vec{W}-\vec{w}))_i(x_j,t_k)| \leq C[M^{-1}+N^{-1}].
\end{equation}
\end{lemma}
\begin{lemma}\label{5pest3} Let the assumptions \eqref{5p3} - \eqref{5p5} hold.
Then, on $\Omega^{N,M}$, for each $1 \leq i \leq n$, the following
estimates hold
\begin{equation*}\label{5p1819} |(L^{N,M}(\vec{W}-\vec {w}))_i(x_j,t_k)|\leq
C(M^{-1}+B_n(x_{j-1})).\end{equation*} 
\end{lemma}
\noindent The following theorem gives the estimate of the error in the singular component $\vec W$.
\begin{theorem} \label{5psingularerrorthm} Let the assumptions \eqref{5p3} - \eqref{5p5} hold. Let $\vec w$ denote the singular component of the solution of 
the problem \eqref{5p1}, \eqref{5p1aa} and $\vec W$ be the singular component of the
 solution of the problem \eqref{5p161}, \eqref{5p161a}.  Then
\begin{equation*}\;\;  \parallel \vec{W}-\vec{w} \parallel  \leq C(M^{-1}+N^{-1}\ln N). \end{equation*}
\end{theorem}
\noindent\textbf {Proof:} 
From the expression \eqref{3175g},
\begin{equation}\label{5p179bb} 
\begin{array}{rcl} 
|(\beta_0^{N,M}(\vec{W}-\vec w))_i(0,t_k)| &\leq&C\sqrt{\eps_i}(x_{1}-x_{0})\ds\max_{s\;\in\;[x_{0},\;x_{1}]}|\dfrac{\partial^2 w_i}{\partial
x^2}(s, t_k)|\\
&\leq&C N^{-1}\ln N,
 \end{array}
 \end{equation}
From the expression \eqref{3175ag},
\begin{equation}\label{5p179bc} 
\begin{array}{rcl} 
|(\beta_1^{N,M}(\vec{W}-\vec w))_i(1,t_k)| &\leq&C\sqrt{\eps_i}(x_{N}-x_{N-1})\ds\max_{s\;\in\;[x_{N-1},\;x_{N}]}|\dfrac{\partial^2 w_i}{\partial
x^2}(s, t_k)|\\
&\leq&C N^{-1}\ln N.
 \end{array}
 \end{equation}
Thus from \eqref{5p179bb}, \eqref{5p179bc} and \eqref{5p184a},
\begin{equation}\label{5p179b} 
\begin{array}{c}
|(\beta_0^{N,M}(\vec{W}-\vec w))_i(0,t_k)| \leq (\beta_0^{N,M}\vec\Phi)_i(0,t_k),\;\;
|(\beta_1^{N,M}(\vec{W}-\vec w))_i(1,t_k)| \leq (\beta_1^{N,M}\vec\Phi)_i(1,t_k),\\
|(\vec{W}-\vec w)_i(x_j,0)| \leq \Phi_i(x_j,0).
 \end{array}
 \end{equation}
In the remaining portion, it is shown that for all $i,j,k$, 
\begin{equation}\label{5p1821} |(L^{N,M}(\vec{W}-\vec{w}))_i (x_j,t_k)|\leq (L^{N,M}\vec{\Phi})_i
(x_j,t_k).
\end{equation} This is proved for each mesh point $x_j \in \Omega^N_x$ by considering separately eight subintervals
\begin{equation*}
 \begin{array}{lcl}
(a)\; (0,\sigma_1),&& \;\;\;(e)\; [1/2,1-\sigma_n],\\
(b)\; [\sigma_1,\sigma_2), &&\;\;\; (f)\; (1-\sigma_{m+1},1-\sigma_m],\; 2 \leq m \leq n-1,\\
(c)\; [\sigma_m,\sigma_{m+1}),\; 2 \leq m \leq n-1,&&\; \;\;\;(g)\; (1-\sigma_2,1-\sigma_{1}]\\
(d)\; [\sigma_n,1/2),&& \;\;\;(h)\; (1-\sigma_1,1).
 \end{array}
\end{equation*}
(a) $x_j\in(0,\sigma_1):\;$Clearly $x_j \notin J$ and $x_{j+1} - x_{j-1} \leq C\sqrt{\eps_1}N^{-1}\ln N.$\\ 
Then, Lemma \ref{5pwwe} and expression \eqref{5p184} give \eqref{5p1821}.\\ Similar arguments hold for the case (h).\\\\
(b) $x_j\in[\sigma_1,\sigma_2):\;$There are 2 possibilities: (b1) $d_1=0$ and (b2) $d_1>0.$ \\\\
(b1) Since $\sigma_1=\dfrac{\sigma_2}{2} $ and the mesh is uniform in $(0,\sigma_2)$ it follows that $x_j \notin J,$ and $x_{j+1} - x_{j-1} \leq C\sqrt{\eps_1}N^{-1}\ln N.$ Then Lemma \ref{5pwwe} and expression \eqref{5p184} give \eqref{5p1821}.\\
(b2)  Either $x_j \notin J$ or $x_j \in J.$ \\
If $x_j \notin J$ then $x_{j+1} - x_{j-1} \leq C\sqrt{\eps_2}N^{-1}\ln N$ and by Lemma \ref{5ps1} $B_1(x_{j-1})\leq B^L_1(x_{j-1})\leq B^L_1(\sigma_2-h^-_2) \leq B^L_1(\sigma_1-h^-_1) \leq CN^{-2},$ so \eqref{5p1815} of Lemma \ref{5pdwd1} with $r=1$ and expression \eqref{5p184} give \eqref{5p1821}.\\
On the other hand, if $x_j \in J$, the expression \eqref{5p1816} of Lemma \ref{5pdwd1} with $r=1$ and expression \eqref{5p185} give \eqref{5p1821}.\\ Similar arguments hold for the case (g).\\\\
(c) $x_j\in[\sigma_m,\sigma_{m+1}):\;$There are 3 possibilities: \\
(c1) $d_1=d_2=\dots =d_m=0,$ \\
(c2) $d_r>0$ and $ d_{r+1}=\;\dots\;=d_m=0$ for some $r,\; 1 \leq r \leq m-1$ and\\
(c3) $d_m>0$.\\\\
(c1) Since the mesh is uniform in $(0,\sigma_{m+1})$, it follows that $x_j \notin J$ and $x_{j+1} - x_{j-1} \leq C\sqrt{\eps_1}N^{-1}\ln N.$  Then Lemma \ref{5pwwe} and expression \eqref{5p184} give \eqref{5p1821}.\\
(c2) Either $x_j \notin J$ or $x_j \in J.$ \\
If $x_j \notin J$ then $\sigma_{r+1}=C\sigma_{m+1},\; \;x_{j+1} - x_{j-1} \leq C\sqrt{\eps_{m+1 }}N^{-1}\ln N$ and by Lemma \ref{5ps1} $B_r(x_{j-1})\leq B^L_r(x_{j-1}) \leq B^L_r(\sigma_m-h^-_m)\leq B^L_r(\sigma_r-h^-_r) \leq CN^{-2}.$  Thus expression \eqref{5p1815} of Lemma \ref{5pdwd1}  and expression \eqref{5p184} give \eqref{5p1821}.\\
On the other hand, if $x_j \in J$, then $x_j=\sigma_m.$ The expression \eqref{5p1816} of Lemma \ref{5pdwd1} with $r=m$ and expression \eqref{5p185} give \eqref{5p1821}.\\
(c3) Either $x_j \notin J$ or $x_j \in J.$ \\
If $x_j \notin J$ then $x_{j+1} - x_{j-1} \leq C\sqrt{\eps_{m+1}}N^{-1}\ln N.$ From \ref{5ps1} $B_m(x_{j-1}) \leq B^L_m(x_{j-1})\leq B^L_m(\sigma_m-h^-_m) \leq CN^{-2}.$
Expression \eqref{5p1815} of Lemma \ref{5pdwd1} with $r=m$ and expression \eqref{5p184} give \eqref{5p1821}.\\
On the other hand, if $x_j=\sigma_m.$ Expression \eqref{5p1816} of Lemma \ref{5pdwd1} with $r=m$ and expression \eqref{5p185} give \eqref{5p1821}.\\ Similar arguments hold for the case (f).\\\\
(d) $x_j\in[\sigma_n,1/2):\;$There are 3 possibilities: \\
(d1) $d_1=\;\;\dots \;\;=d_n=0,$\\ 
(d2) $d_r>0$ and $ d_{r+1}=\;\dots\;=d_n=0$ for some $r,\; 1 \leq r \leq n-1$ and\\
(d3) $d_n>0$.\\\\
(d1) Since the mesh is uniform in $\Omega^N_x$, it follows that $x_j \notin J,\;x_{j+1} - x_{j-1} \leq C\sqrt{\eps_1}N^{-1}.$ Then Lemma \ref{5pwwe} and expression \eqref{5p184} give \eqref{5p1821}.\\
(d2) Either $x_j \notin J$ or $x_j \in J.$ \\
If $x_j \notin J$ then $x_{j+1} - x_{j-1} \leq C\sqrt{\eps_{r+1}}N^{-1}$ and by Lemma \ref{5ps1}, $B_r(x_{j-1}) \leq B^L_r(x_{j-1})\leq B^L_r(\sigma_n-h^-_n)\leq B^L_r(\sigma_r-h^-_r) \leq CN^{-2}.$  The expression \eqref{5p1815} of Lemma \ref{5pdwd1} and expression \eqref{5p184} give \eqref{5p1821}.\\
On the other hand, if $x_j \in J$, then $x_j = \sigma_n.$\\
The expression \eqref{5p1816} of Lemma \ref{5pdwd1} and expression \eqref{5p185} give \eqref{5p1821}.\\
(d3) By Lemma \ref{5ps1} with $r=n$, $B_n(x_{j-1}) \leq B^L_n(x_{j-1})\leq B^L_n(\sigma_n-h^-_n) \leq CN^{-2}.$ \\
Then Lemma \ref{5pest3} and expression \eqref{5p184} give \eqref{5p1821}.\\ Similar arguments hold for the case (e).\\
By using comparision principle, the required result is established from \eqref{5p179b} and \eqref{5p1821}.\\

\noindent The following theorem gives a parameter uniform bound which is first order in time and essentially
first order in space for the convergence of the discrete solution.
\begin{theorem}\label{5pmain} Let the assumptions \eqref{5p3} - \eqref{5p5} hold. Let $\vec u$ denote the solution of the problem \eqref{5p1}, \eqref{5p1aa}
and $\vec U$ denote the solution of the problem \eqref{5p161}, \eqref{5p161a}.  Then
\begin{equation*}\;\;  \parallel \vec{U}-\vec{u} \parallel  \leq C(M^{-1}+N^{-1}\ln N). \end{equation*}
\end{theorem}
\noindent \textbf{Proof.}
An application of the triangular inequality and the results of
Theorem \ref{5psmootherrorthm} and Theorem \ref{5psingularerrorthm} lead to the required result.
\section{Numerical Illustration}\label{c5s38}
\qquad The numerical method proposed above is illustrated through the example presented in this section. 
The method proposed above is applied to solve the
problem and the parameter-uniform order of convergence and the parameter-uniform error constants are computed. 
To get the order of convergence in the variable $t$ seperately, a Shishkin mesh is considered for $x$ and 
the resulting problem is solved for various uniform meshes with respect to $t$. In order to get the order of convergence in the variable $x$ seperately, a uniform mesh
is considered for $t$ and the resulting problem is solved for various piecewise uniform Shishkin meshes with respect to
$x$. The two-mesh algorithm  for a vector problem which is a variant of the one found in \cite{12} is applied to get parameter-uniform order of 
convergence and the error constants. The numerical results are presented in Tables \ref{t5p1} and Table \ref{t5p2}.\\

{\noindent\bf Example 1 }
Consider the problem 
\[\frac{\partial \vec u}{\partial t} -E\frac{\partial^2 \vec u}{\partial x^2} + A \vec u= \;\vec{f}\text{ on }
(0,1) \times (0,1],\]\[(\vec u-E_*\dfrac{\partial \vec u}{\partial x})(0,t)=\vec{\phi}_L,\;\;
(\vec u+E_*\dfrac{\partial \vec u}{\partial x})(1,t)=\vec{\phi}_R,\;\;\vec u(x,0)= \vec{\phi}_B\]
where \;$ E = diag(\eps_1, \eps_2)$,\; $A=$
$\begin{pmatrix}
4+3t\;\;& -1\\
-1\;\;& 4+3t
\end{pmatrix},\;
\vec f = 
\begin{pmatrix}
2+e^{3t}\\
2+e^{3t}
\end{pmatrix},\;
\vec{\phi}_L = 
 \begin{pmatrix}1+t^8\\1+t^8
\end{pmatrix},\;
\vec{\phi}_R = 
\begin{pmatrix}1+t^8\\1+t^8\end{pmatrix}, \;
\vec{\phi}_B = \begin{pmatrix}1\\1\end{pmatrix}.$\\\\
For various values of $\eps_1$ and $\eps_2$, the maximum errors, the $\vec\eps$- uniform order of convergence and the $\vec\eps$-uniform error constant are computed.
Fixing a Shishkin mesh on $[0,1]$ with $128$ points horizontally, the problem is solved by the method 
suggested above. The order of convergence and the error constant for $\vec u$ are calculated for $t$ using two-mesh algorithm and the
results are presented in Table \ref{t5p1}. A uniform mesh on $[0,1]$ with $32$ points vertically is considered and the order of
convergence and the error constant for $\vec u$ in the variable $x$ using two-mesh algorithm are calculated and the results are presented in Table \ref{t5p2}.

\begin{figure} 
\begin{minipage}{0.45\textwidth}\caption{\label{f5p1}}
\begin{center}The numerical approximation of $\vec u$ \\for $\eps_1=2^{-15},\;\eps_2=2^{-14}$ and $M=32$\end{center} 
\includegraphics[width=\textwidth]{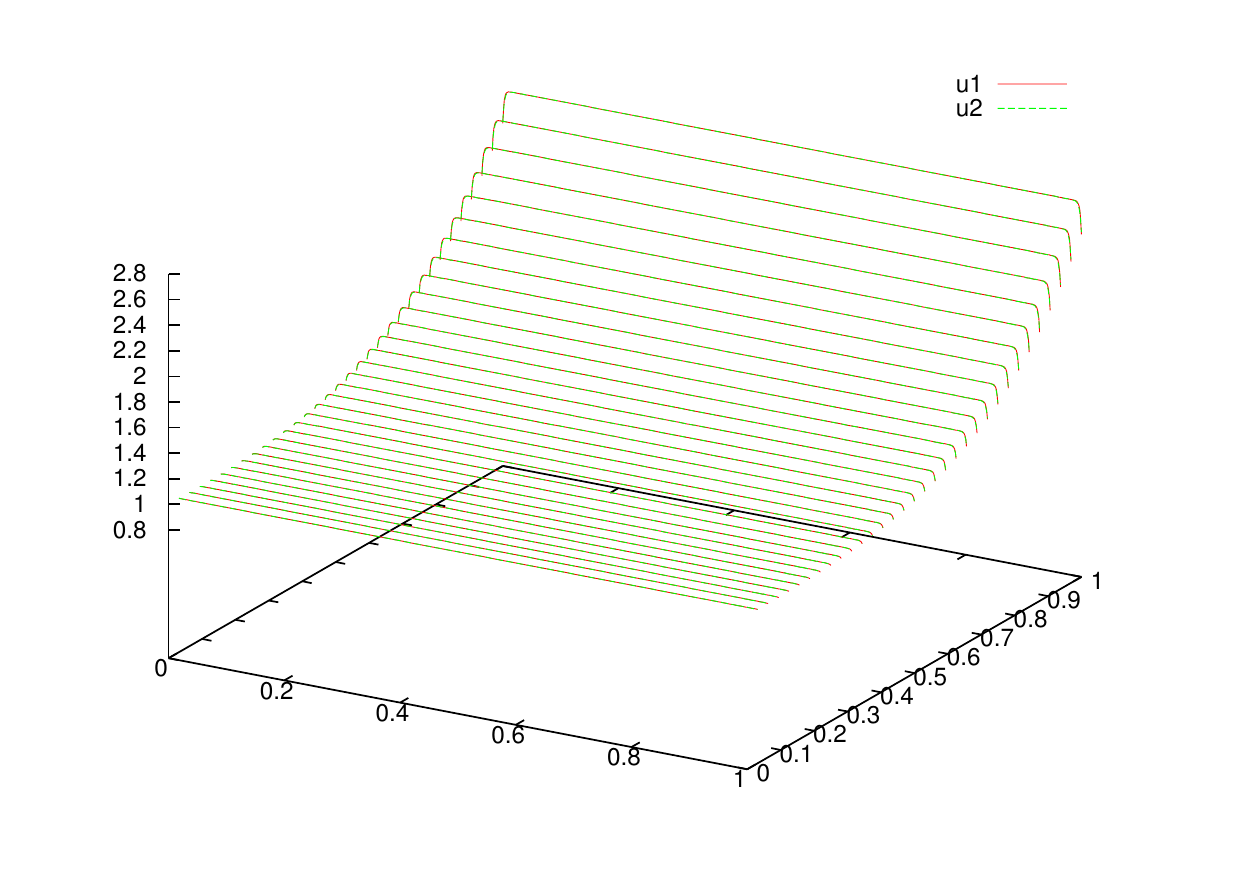} 
\end{minipage}\;\;
\begin{minipage}{0.45\textwidth}\caption{\label{f5p2}}
\begin{center}The numerical approximation of $ \vec u$ \\for $\eps_1=2^{-15},\;\eps_2=2^{-14}$ and $N=128$\end{center} 
\includegraphics[width=\textwidth]{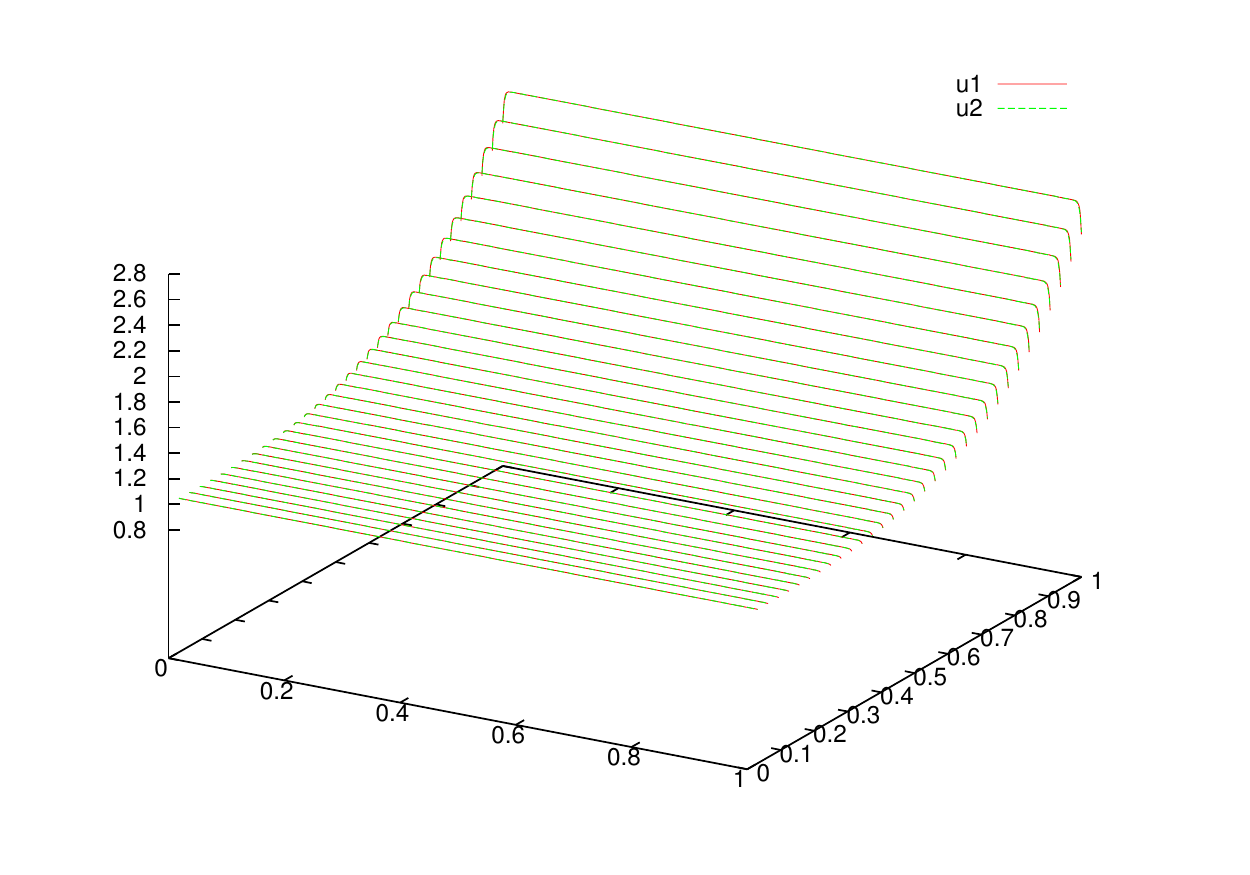} 
\end{minipage} 
\end{figure}

It is evident from the Figures \ref{f5p1} and \ref{f5p2} that the solution $\vec u$ exhibits
parabolic twin boundary layers at $ (0,t) $ and $ (1,t), \; 0 \leq t \leq 1. $ Further, the $t$- order of convergence and the $x$- order of convergence of
the numerical method presented in Table \ref{t5p1} and Table \ref{t5p2} agree with the theoretical result. 

\begin{table}[!ht]
\begin{center}
 \caption{\label{t5p1}Values of \;\;$D_{\vec{\varepsilon}}^N, D^N, p^{N}, p^* \text{ and }C_{p^*}^N\text{ for }\;\varepsilon_1=\dfrac{\eta}{16},\;\varepsilon_2=\dfrac{\eta}{8},\; \alpha = 2.9,\;\text{and} \;N=128.$}
\end{center}	
\begin{center}
\begin{tabular}{|l|l|l|l|l|}
\hline         
\multicolumn{1}{|c|}{$\eta$} &\multicolumn{4}{|c|}{Number of mesh points $M$}\\
\cline{2-5}& \;\;\;32 &\;\;\;64&\;\;\;128  & \;\;\; 256 \\\hline
 \;\;$2^{-7}$  &  0.153E-01   &    0.783E-02  &     0.397E-02   &    0.199E-02\\\hline
  \;\;$2^{-8}$  &  0.155E-01   &    0.788E-02   &    0.397E-02   &    0.200E-02\\\hline
  \;\;$2^{-9}$  &  0.156E-01   &    0.790E-02   &    0.398E-02   &    0.200E-02\\\hline
  \;\;$2^{-10}$  &   0.156E-01  &     0.790E-02  &     0.398E-02  &     0.200E-02\\\hline
  \;\;$2^{-11}$  &  0.156E-01   &    0.790E-02   &    0.398E-02   &    0.200E-02\\\hline\hline
 
  \; $D^N$  &    0.156E-01   &    0.790E-02    &   0.398E-02   &    0.200E-02\\\hline
     \; $p^N$ &  0.980E+00   &    0.990E+00    &   0.995E+00&\\\hline
    \; $C_{p^*}^N$ &  0.945E+00   &    0.945E+00  &     0.939E+00    &   0.929E+00\\\hline\hline

\multicolumn{5}{|c|}{Computed $t$-order of $\vec{\varepsilon}-$uniform convergence, $p^*$ = 0.9803767} \\\hline
\multicolumn{5}{|c|}{Computed $\vec{\varepsilon}-$uniform error constant, $C_{p^*}^*$ = 0.9451866 }\\\hline
\end{tabular}
\end{center}
\end{table}
\begin{table}[!ht]
\begin{center}
 \caption{\label{t5p2}Values of \;\;$D_{\vec{\varepsilon}}^N, D^N, p^{N}, p^* \text{ and }C_{p^*}^N\text{ for }\;\varepsilon_1=\dfrac{\eta}{16},\;\varepsilon_2=\dfrac{\eta}{8},\; \alpha = 2.9,\;\text{and} \;M=32.$ }
\end{center}
\begin{center}
\begin{tabular}{|l|l|l|l|l|}
\hline         
\multicolumn{1}{|c|}{$\eta$} &\multicolumn{4}{|c|}{Number of mesh points $N$}\\
\cline{2-5}& \;\;\;32 &\;\;\;64&\;\;\;128  & \;\;\; 256 \\\hline
 \;\;$2^{-7}$    &   0.530E-01   &    0.339E-01  &     0.172E-01   &    0.689E-02\\\hline
    \;\;$2^{-8}$  &   0.530E-01  &     0.339E-01  &     0.172E-01  &     0.689E-02\\\hline
   \;\;$2^{-9}$   &  0.530E-01   &    0.339E-01   &    0.172E-01   &    0.689E-02\\\hline
  \;\;$2^{-10}$   &  0.530E-01   &    0.339E-01   &    0.172E-01   &    0.689E-02\\\hline
  \;\;$2^{-11}$   &  0.530E-01   &    0.339E-01   &    0.172E-01   &    0.689E-02\\\hline\hline
 
  \; $D^N$  &     0.530E-01    &   0.339E-01  &     0.172E-01   &    0.689E-02\\\hline
     \; $p^N$ &    0.644E+00   &    0.977E+00  &     0.132E+01&\\\hline
    \; $C_{p^*}^N$ &  0.137E+01  &     0.137E+01  &     0.109E+01   &    0.679E+00\\\hline\hline

\multicolumn{5}{|c|}{Computed $x$-order of $\vec{\varepsilon}-$uniform convergence, $p^*$ = 0.6436486} \\\hline
\multicolumn{5}{|c|}{Computed $\vec{\varepsilon}-$uniform error constant, $C_{p^*}^*$ = 1.371360 }\\\hline
\end{tabular}
\end{center}
\end{table}
 
\section*{Acknowledgment}

The first author sincerely thanks the University Grants Commission, New Delhi, India, for the financial support through the Rajiv Gandhi National
Fellowship to carry out this work.


\end{document}